\documentclass[reqno,11pt, letterpaper]{amsart}
\PassOptionsToPackage{fleqn}{amsmath}
\RequirePackage{amsthm}
\RequirePackage{amsmath}
\RequirePackage{latexsym}
\RequirePackage{amssymb}
\RequirePackage{amscd}
\RequirePackage{epsfig}
\RequirePackage{graphics}
\RequirePackage{ifthen}
\RequirePackage{varioref}
\usepackage{changes}
\usepackage{soul}
\usepackage{xcolor}
\usepackage{color}
\usepackage[misc]{ifsym}
\usepackage{bm}
\usepackage{lipsum}
\usepackage{subfigure}

\usepackage{epstopdf}

\definecolor{darkgreen}{rgb}{0.0625,0.64,0.0625}
\ifpdf
  \usepackage[
    pdftex,
    colorlinks,%
    linkcolor=blue,citecolor=blue,urlcolor=blue,
    hyperindex,%
    plainpages=false,%
    bookmarksopen,%
    bookmarksnumbered%
  ]{hyperref}
\else
  \usepackage{hyperref}
\fi
\usepackage{bookmark}
\allowdisplaybreaks[1]

\def\rd{\mathrm{d}}

\def\R{{\mathbb R}}

\theoremstyle{plain}
\newtheorem{thm}{Theorem}[section]
\newtheorem{lem}[thm]{Lemma}
\newtheorem{prop}[thm]{Proposition}

\theoremstyle{definition}
\newtheorem{rem}[thm]{Remark}

\newtheorem{defn}[thm]{Definition}

\def\ro#1{{\rm #1}}

\def\Bbb#1{{\mathbb#1}}

\def\R{\Bbb R}

\def\Rx{\R\mkern1mu^}

\def\Rn{\Rx n}

\def\AG#1{\ro{A}(#1)}  
\def\AGR#1{\AG{#1,\R}} 
\def\SA#1{\ro{SA}(#1)} 
\def\SAR#1{\SA{#1,\R}} 
\def\GL#1{\ro{GL}(#1)} 
\def\GLR#1{\GL{#1,\R}} 
\def\SL#1{\ro{SL}(#1)} 
\def\SLR#1{\SL{#1,\R}}  

\def\semidirect{\ltimes}

\numberwithin{equation}{section}
\newcommand\blfootnote[1]{%
 \begingroup
 \renewcommand\thefootnote{}\footnote{#1}%
 \addtocounter{footnote}{-1}%
 \endgroup
 }

\begin{document}

\title[Equi-centro-affine extremal hypersurfaces in ellipsoid]{Equi-centro-affine extremal\\ hypersurfaces in ellipsoid}

    \author[Y. Yang]{Yun Yang}
    \address{Yun Yang\newline\indent
     Department of Mathematics, Northeastern University, Shenyang, 110819, P.R. China}
    \email[Corresponding author]{yangyun@mail.neu.edu.cn}

    \author[C. Z. Qu]{Changzheng Qu${}^*$}\blfootnote{${}^*$~Corresponding author: quchangzheng@nbu.edu.cn}
    \address{Changzheng Qu\newline\indent
     School of Mathematics and Statistics, Ningbo University, Ningbo, 315211, P.R. China}
    \email{quchangzheng@nbu.edu.cn}

\begin{abstract}
    This paper explores equi-centro-affine extremal hypersurfaces in an ellipsoid.
    By analyzing the evolution of invariant submanifold flows under centro-affine unimodular transformations, we derive the first and second variational formulas for the associated  invariant area.
    Stability analysis reveals that the circles with radius $r=\sqrt{6}/3$ on $\mathbb{S}^2(1)$ are characterized as being equi-centro-affine maximal.
    Furthermore, we provide a detailed classification of the compact isoparametric equi-centro-affine extremal hypersurfaces on $(n+1)$-dimensional sphere, as well as the generalized closed equi-centro-affine extremal curves on $2$-dimensional sphere.
    These curves are shown to belong to a family of transcendental curves $\mathrm{x}_{p,q}$  ($p,q$ are two coprime positive integers satisfying that $1/2<p/q<1$ ).
    Additionally, we establish an equi-centro-affine version of isoperimetric inequality ${}^{ec}\hspace{-1mm}L^3\leq (4\pi-A)(2\pi-A)A$ on $\mathbb{S}^2(1)$.
\end{abstract}

\subjclass[2010]{53A15, 53A55, 53A10, 58E12.}

\keywords{equi-centro-affine geometry;\; Euler-Lagrange equation;\; variational formula;\;  isoperimetric inequality}


\maketitle
\section{Introduction}
In general,  minimal submanifolds serve as the critical points of the area or volume functional \cite{bla0,bla1,cm}, which also emerge as static solutions of the mean curvature flow \cite{hui-84,zhu}.
On one hand, minimal submanifolds, occupying a prominent position in global differential geometry, have been studied extensively, yielding a plethora of fascinating results (refer to  \cite{bre12,bre131, bre21,bre} for novel and significant achievements).
On the other hand, various delicate methods such as variational method, Min-Max theory, have been proposed to construct minimal submanifolds \cite{mn14,mn16,mn17},
and subsequently, the analysis of the second variation can provide valuable insights into the stability, as detailed in \cite{cl0,cl}.
Certainly, the theory of minimal submanifolds, being a pivotal theme in geometric analysis, has been affirmed as a powerful and essential instrument in mathematics \cite{bre13}.
Sustained development addressing  innovative and substantial accomplishments in minimal submanifolds, besides being important in its own right, may also greatly
enhance potential applications in related disciplines such as computer vision, probability and general relativity.
In this paper, we focus on the equi-centro-affine extremal hypersurfaces $M^n$ in an ellipsoid $N^{n+1}\subset\R^{n+2}$.

The problem of finding a minimal surface with a prescribed boundary has engaged such prolific mathematicians as Lagrange and Euler, as well as, Plateau. It was Plateau who first delved into the study of the surface obtained in the form of a soap film stretched on a wire framework, a physical example of a minimal surface, and this problem has become known as the famous Plateau's problem \cite{cour, dou, rad}. Other noteworthy examples of minimal surfaces include helicoids and catenoids \cite{cm0}. The Bernstein problem also holds a pivotal role in the theory of minimal submanifolds \cite{alm,bgg}, and numerous investigations have emerged focusing on Bernstein type problems in Euclidean spaces (see \cite{fs,ssy,sol} and the references therein).

A natural generalization is to study minimal
surfaces in Riemannian manifolds, deviating from the conventional setting of $\R^{n}$,
with an interesting instance being the $n$-dimensional sphere, denoted as $\mathbb{S}^n$ \cite{bre12, bre131,bre13,cck}.
Note that a crucial distinction from the ${\mathbb R}^n$ scenario lies in the fact that:
every minimal submanifold in $\R^n$ inherently possesses a non-compact nature, whereas in $\mathbb{S}^n$, one can encounter closed minimal submanifolds.

Let $M^n$ be an $n$-dimensional Riemannian manifold that is immersed isometrically into a space form $N^{n+1}$ with constant curvature $c$.
The immersion's principal curvatures are denoted by $\kappa_1,\cdots,\kappa_n$. For $r$ ranging from $0$ to $n$, let $S_r$
represent the $r$-th elementary symmetric polynomial, which is given by the sum of all possible products of $r$ distinct principal curvatures:
$\displaystyle \sum_{i_1<i_2<\cdots<i_r}\kappa_{i_1}\kappa_{i_2}\cdots\kappa_{i_r}$.
Reilly \cite{rei} extended variational problems concerning the area or volume functional to any
smooth function $f(S_1,\cdots,S_n)$ defined on the manifold $M^n$, specifically integrals of the form  $\displaystyle \int_Mf(S_1,\cdots,S_n)dV$, where $dV$ represents the volume element on $M^n$.
The Willmore submanifold \cite{bry, lhz}  is an extremal submanifold in terms of the Willmore
functional $\displaystyle \int_M (S-nH^2)^{n/2}dV$, where $S$  denotes the square of the length of the second fundamental
form, $H$ signifies the mean curvature of $M^n$. Notably, it remains invariant under M\"obius (or
conformal) transformations of $\mathbb{S}^{n+1}(1)$.

The variational issue in the affine setting is a little more complicated (refer to \cite{wan1,wan2} for the variational problems respect to the affine arc length).
The exploration of affine differential geometry is grounded in the Lie group $\AGR n = \GLR n \semidirect \Rn$  which includes affine transformations of the form $x \longmapsto Ax+b$, $A\in \GLR n$, $b\in \Rn$ acting on $x \in \Rn$
(refer to Nomizu and Sasaki \cite{ns} and Simon \cite{sim} for details). Analogously, equi-affine geometry is confined to the subgroup $\SAR n = \SLR n \semidirect \Rn$ of volume-preserving affine transformations. Centro-affine differential geometry refers to the subgroup of the affine transformation group that keeps the origin fixed, which is closely related to the geometry induced by the general linear group $x \longmapsto Ax$, $A\in \GLR n$, $x \in \Rn$.
Furthermore, equi-centro-affine differential geometry arises in connection with the subgroup $\SLR n$ of volume-preserving linear transformations.

In equi-affine differential geometry, the ambient
space $\R^{n+2}$ has a flat affine connection $D$ and the usual determinant function
 is regarded as a parallel volume element.
Let $\mathbf{x}: N^{n+1}\rightarrow \R^{n+2}$ be a local embedding of a smooth hypersurface,
 and $\xi$ be the affine normal field to $\mathbf{x}$. The equi-affine structure equations of $\mathbf{x}$ may be written as (see \cite{lt,ns} for more details)
    \begin{align*}
      \mathbf{x}_{ij}&=\bar{g}_{ij}\xi+(\bar{\Gamma}^k_{ij}+C^k_{ij})\mathbf{x}_k,\\
      \xi_i&=-A^k_i\mathbf{x}_k,
    \end{align*}
    where $\bar{g}_{ij}$ is the equi-affine metric, $\bar{\Gamma}^k_{ij}$ are the Christoffel  symbols of the metric $\bar{g}_{ij}$, $C^k_{ij}$ is called cubic form, and $A^k_i$ is the equi-affine shape operator.

The Euclidean inner product $\langle\cdot,\cdot\rangle$ on $\R^{n+2}$ induces the metric ${}^{e}\!g_{ij}$ and the Euclidean second fundamental
form ${}^{e}\!h_{ij}$ on $\mathbf{x}$.
According to the literature \cite{lt}, the connection between the equi-affine and Euclidean metrics is given by
\begin{equation}\label{equi-metric}
         \bar{g}_{ij}=\frac{{}^{e}\!h_{ij}}{\phi},
\end{equation}
where $\displaystyle \phi=\left(\frac{\det{{}^{e}\!h_{ij}}}{\det{{}^{e}\!g_{ij}}}\right)^{1/(n+2)}$.

     The extremal submanifolds in equi-affine space $\R^{n+2}$ are the critical points of the equi-affine invariant area functional given by $\displaystyle \int_N\sqrt{\left|\det\bar{g}_{ij}\right|}d\mu_N$.
     Note that, commonly in the literature, the term ``affine geometry'' is used interchangeably with ``equi-affine geometry''.

   In approximately 1977, Chern \cite{che} conjectured that an affine maximal graph of a smooth, locally uniformly convex function on two-dimensional Euclidean space, $\R^2$, is necessarily a paraboloid.
   The two-dimensional Chern's conjecture was fully resolved by Trudinger and Wang in their celebrated paper \cite{tw-1}. However, the higher-dimensional Bernstein problem  remains unsolved. Later, Li and Jia \cite{lj}, and also Trudinger and Wang \cite{tw-11}, proved Calabi's conjecture for two-dimensions separately, using distinctively different approaches. Additionally, Trudinger and Wang investigated the Plateau problem for affine maximal hypersurfaces,
   which serves as the analogous affine-invariant counterpart of the classical Plateau problem for minimal
   surfaces \cite{tw-2}. In \cite{wan}, Wang declared that the concept of an affine maximal surface in affine geometry  mirrors  that of minimal surface in Euclidean geometry (Calabi \cite{cal} advocated for the terminology ``affine maximal" as the
   second variation of the affine area functional is negative).
   The affine Bernstein problem and the affine Plateau problem, as introduced in \cite{cal,cal-1,che}, constitute two fundamental issues of affine maximal submanifolds.

   In this paper, we explore the equi-centro-affine maximal hypersurfaces $M^n$ immersed in the ellipsoid $N^{n+1}$ which is centered at the origin in $\R^{n+2}$.
   Certainly, it is universally acknowledged that an ellipsoid in the equi-affine setting is an affine sphere with vanishing cubic form.
   Under an equi-affine transformation, the ellipsoid $N^{n+1}$ may be formulated as
    \begin{align}
       \mathbf{x}=&R\Big\{\cos r,\sin r\cos\theta^1,\;\cdots,\; \sin r\sin\theta^1 \cdots\sin\theta^{n-1}\cos\theta^{n},\nonumber\\
                  &\qquad\qquad\qquad\qquad\qquad\qquad\qquad\qquad\qquad \sin r\sin\theta^1\cdots\sin\theta^{n}\Big\}.\label{Ellipsoid}
    \end{align}
    where $R=\left(a_0a_1\cdots a_{n+1}\right)^{1/(n+2)}$, and $a_0,a_1,\cdots, a_{n+1}$ represent the lengths of the semi-axes of the ellipsoid. Denote $\varrho=R^{\frac{n+2}{n+3}}$. Then by \eqref{equi-metric}, the equi-affine metric of the ellipsoid $\mathbf{x}$ can be written as
    \begin{equation}\label{ECA}
        \bar{g}=\varrho^2
        \left(
          \begin{array}{ccccc}
            1 & 0  & \cdots & 0 \\
            0 & \sin^2r &  \cdots & 0 \\
            0 & 0  & \cdots & 0 \\
            \cdots & \cdots & \cdots & \cdots \\
            0 & 0  & \cdots & \sin^2r\sin^2\theta^1\cdots\sin^2\theta^{n-1} \\
          \end{array}
        \right).
    \end{equation}
    Since $R$ and $\varrho$ are equal up to a constant scaling factor, it is permissible to interchange $\varrho$ with $R$ in equi-affine metric \eqref{ECA} of the ellipsoid $N^{n+1}$.

    Alternatively,  $M^n$ is a submanifold of codimension two in $\R^{n+2}$, endowed with an equi-centro-affine structure, as detailed in \cite{ns93,ns,wal,yl}.
    Then for any local oriented basis $\sigma =
\{E_{1},E_{2},\cdots,E_{n}\}$ of $TM$ with dual
basis $\{\theta^{1},\theta^{2},\cdots,\theta^{n}\}$, we introduce the term
\begin{eqnarray}\label{Gdef}
G:=[E_{1}(x),\cdots,E_{n}(x),x,\rd^2
x]=G_{ij}\theta^{i}\otimes\theta^{j},
\end{eqnarray}
under the assumption that $G$ is nondegenerate,
where the bracket notation $[\cdots]$ is employed to denote the standard determinant in $\R^{n+2}$, and
  $$G_{ij}:=[E_{1}(x),\cdots,E_{n}(x),x,E_{i}E_{j}(x)]$$ is a symmetric $2$-form.
Moreover, we may verify that
\begin{equation}\label{Metricdef}
\tilde{g}:=\tilde{g}_{ij}\theta^{i}\otimes\theta^{j},\quad\quad\displaystyle \tilde{g}_{ij}=|\det (G_{pq})|^{-\frac{1}{n+2}}G_{ij}
\end{equation}
 is independent of the choice of the basis $\sigma$ and thus
a globally defined symmetric $2$-form, which also is invariant up to the equi-centro-affine transformations in $\R^{n+2}$.
According to \cite{yl}, $\{\tilde{g}_{ij}\}$ is identified  as  an equi-centro-affine metric  associated with the immersion $\mathbf{x}:M^n\rightarrow \R^{n+2}$.
 Let $\Delta_{\tilde{g}}$ denote the Laplacian of $\tilde{g}$. $\displaystyle \{\mathbf{x},\frac{\Delta_{\tilde{g}} \mathbf{x}}{n}\}$ is characterized as the equi-centro-affine normalization
of $\mathbf{x}:M^n\rightarrow \R^{n+2}$.

The focus of the present question shifts to elucidating the extremal equi-centro-affine hypersurfaces of the ellipsoid $N^{n+1}$ endowed with the equi-affine metric in $\R^{n+2}$,  which arise from the application of the equi-centro-affine metric defined on $M^n$.

Given that $M^n$ is furthermore a submanifold residing in the ellipsoid $N^{n+1}$ centered  at origin of $\R^{n+2}$. The equi-centro-affine invariant area function $\displaystyle \int_M\sqrt{\det\tilde{g}_{ij}}d\mu_M$  can be denoted by (see Section \ref{sec-eca})
\begin{align*}
       \Sigma_{eca}=\varrho^{\frac{n}{n+2}}\int_M S_n^{\frac{1}{n+2}}\sqrt{\det(g_{ij})}d\mu_M,
\end{align*}
where $g_{ij}$ is the metric of $M^n$ induced from $\bar{g}_{ij}$.
Note that, as a consequence of the nondegeneracy condition of the metric $\tilde{g}$,  $S_n\neq0$ holds true on the submanifold $M^n$.
Furthermore, when $n$ takes on an even number,   the same nondegeneracy condition ensures that $S_n>0$.
In Section \ref{sec-VF}, we will prove the following theorem:
\begin{thm}\label{thm-fv}
The first variational formula, for the equi-centro-affine invariant area of a hypersurface $M^n$ in the ellipsoid $N^{n+1}$ is
\begin{align*}
       \frac{d}{dt}\Sigma_{eca}=&\varrho^{\frac{n}{n+2}}\int_M\frac{U}{n+2}\Bigg(T_{n-1}{}^{ij}\nabla_{e_i}\nabla_{e_j}\left(S_n^{-\frac{n+1}{n+2}}\right)\\
                           &\qquad +S_n^{-\frac{n+1}{n+2}}\left(\frac{1}{\varrho^2}S_{n-1}-(n+1)S_nS_1\right)\Bigg)dV.
    \end{align*}
\end{thm}
Let us denote
\begin{align}\label{extremal}
    {}^{ec}\hspace{-1mm}H= T_{n-1}{}^{ij}\nabla_{e_i}\nabla_{e_j}\left(S_n^{-\frac{n+1}{n+2}}\right)+S_n^{-\frac{n+1}{n+2}}\left(\frac{1}{\varrho^2}S_{n-1}-(n+1)S_nS_1\right).
  \end{align}
The hypersurface in the ellipsoid is said to be an equi-centro-affine extremal hypersurface if ${}^{ec}\hspace{-1mm}H=0$.
\begin{rem} When $n=2$ and $R=1$, ${}^{ec}\hspace{-1mm}H=0$ reduces to
\begin{align*}
  S_1\left(\Delta_{g}\left(S_2^{-\frac{3}{4}}\right)+S_2^{-\frac{3}{4}}\left(1-3S_2\right)\right)=0,
\end{align*}
where $\Delta_{g}$ denotes the Laplace-Beltrami operator with respect to metric $g$ of $M^n$ induced from $\bar{g}$.
Since $S_2>0$, we conclude that $S_1\neq0$. Thus
\begin{align*}
  \Delta_{g}\left(S_2^{-\frac{3}{4}}\right)+S_2^{-\frac{3}{4}}\left(1-3S_2\right)=0.
\end{align*}
The trivial solution for this equation is $\displaystyle S_2=1/3$. Assume $M^2$ is compact, we obtain a Simon's type equality (see \cite{lhz,simj} for Simon's type integral inequality)
\begin{align*}
  \int_M S_2^{-\frac{3}{4}}\left(1-3S_2\right) dV=0.
\end{align*}
\end{rem}

Let us draw an analogy with the affine maximal hypersurface in $\R^n$. If the affine maximal hypersurface $M^n$ in $\R^{n+1}$ can be represented as the graph of a convex function $u$, then $u$ satisfies
the Monge-Ampere type equation $\displaystyle \Delta_{\bar{g}}\left(\det\left(D^2u\right)^{-1/(n+2)}\right)=0$ and
$\Delta_{\bar{g}}$ denotes the Laplace-Beltrami operator with respect to  equi-affine metric $\bar{g}$. It is noteworthy that the higher-dimensional affine maximal hypersurface problems remain unsolved to date.

The Clifford type hypersurfaces
 \begin{align*}
     M^n=\mathbb{S}^m\left(\sqrt{\frac{m+1}{n+2}}\right)\times\mathbb{S}^{n-m}\left(\sqrt{\frac{n+1-m}{n+2}}\right),\qquad 1\leq m\leq n-1
    \end{align*}
 are examples of equi-centro-affine extremal hypersurfaces in unit sphere (refer to Section \ref{sec-eh} for a comprehensive classification of isoparametric equi-centro-affine extremal hypersurfaces).
 Under the specified conditions $2m=n$ and $m$ being an even number,  the aforementioned Clifford-type  hypersurfaces $M^n$ are Euclidean minimal hypersurfaces and also Willmore hypersurfaces.

In Section \ref{sec-vf2}, we derive the second variational formula for the equi-centro-affine invariant area of a hypersurface $M^n$ immersed in the ellipsoid $N^{n+1}$ (see Theorem  \ref{thm-second}).
Based on this theorem, we embark on an investigation into the stability of curves with respect to the equi-centro-affine arc length when they are considered on the unit sphere $\mathbb{S}^2(1)$. One of our main results is
\begin{thm}\label{thm-cc}
    On the unit sphere $\mathbb{S}^2(1)$, the circle with radius $r=\frac{\sqrt{6}}{3}$ stands out as a unique embedded closed curve that is simultaneously stable and equi-centro-affine maximal.
\end{thm}
\begin{rem}
  According to Theorem \ref{thm-cc},  on $\mathbb{S}^2(1)$,  if there exists a maximizer of the equi-centro-affine arc length among all embedded closed curves, then that curve is necessarily a planar circle with radius  $r=\sqrt{6}/{3}$.
  In this scenario, the equi-centro-affine arc length ${}^{ec}\hspace{-1mm}L$ of any embedded closed curve on $\mathbb{S}^2(1)$ is upper-bounded by ${}^{ec}\hspace{-1mm}L\leq 2^{\frac{4}{3}}\sqrt{3}\pi/{3}$. Let $A$ denote the area enclosed by an embedded closed curve on open hemisphere of $\mathbb{S}^2(1)$.  Assuming further that  a maximizer of the equi-centro-affine arc length exists for curves with a fixed enclosed area, Section \ref{sec-fix} and Theorem \ref{thm-ca} jointly establish a profound equi-centro-affine isoperimetric inequality ${}^{ec}\hspace{-1mm}L^3\leq (4\pi-A)(2\pi-A)A$, with equality achieved solely by planar circles.
\end{rem}
For a closed equi-centro-affine extremal curve on the unit sphere $\mathbb{S}^2(1)$, the progression angle $\Lambda^{\Theta}$ in one period of the curvature is
\begin{align*}
   \Lambda^{\Theta} = \;&\frac{4}{(a+r)\sqrt{(a-c)}}\Pi\left(\frac{a-b}{a+r},\frac{\pi}{2},\frac{a-b}{a-c}\right)\nonumber\\
                    &-\frac{4}{(a-r)\sqrt{(a-c)}}\Pi\left(\frac{a-b}{a-r},\frac{\pi}{2},\frac{a-b}{a-c}\right),
 \end{align*}
where $C_2 > 3\sqrt[3]{4}$, and $\Pi$ is the elliptic integral of the third kind. Additionally, $a, b, c$ are the three distinct solutions of the cubic polynomial $x^3 - C_2x + 4 = 0$, satisfying the specific ordering $a > b > 0 > c$, and $r = \sqrt{C_2}$.
In Appendix B, we will prove that $\Lambda^{\Theta}$ exhibits monotonic decrease with the increase of the parameter $C_2$, and $\pi < \Lambda^{\Theta} < \sqrt{2}\pi$.
In relation to the classification of closed equi-centro-affine extremal curves on the unit sphere $\mathbb{S}^2$, we establish
\begin{thm}
Let $\mathbf{x}$ be a closed equi-centro-affine extremal curve on the unit sphere $\mathbb{S}^2(1)$. Then we have the following possibilities for $\mathbf{x}$:
 \begin{itemize}
   \item[(1)] $\mathbf{x}$ is a planar circle with radius $\sqrt{6}/3$;
   \item[(2)] $\mathbf{x} = \mathbf{x}_{p,q}$ has rotation index $p$ and closes up in $q$ periods of its curvature function.
    The pair $(p, q)$ is not arbitrary and must be such that $p/q$ is defined in the open interval $\left(\frac{1}{2}, \frac{\sqrt{2}}{2}\right)$.
 \end{itemize}
\end{thm}

The rest of this paper is organized as follows. In Section 2, we review the geometry of equi-affine space and equi-centro-affine space, and deduce the equi-centro-affine invariant area function for the hypersurface $M^
n$ immersed in ellipsoid $N^{n+1}$. Section 3 is dedicated to the computation of variational formulas, specifically focusing on the equi-centro-affine area and arc-length of equi-centro-affine curves residing on the ellipsoid. In Section 4, our attention shifts to the investigation of equi-centro-affine extremal hypersurfaces on the unit sphere. Lastly,  in Section 5, the closed equi-centro-affine extremal curves on the sphere are studied. As a conclusion, a classification of closed generalized equi-centro-affine extremal curves on the sphere is presented.
\section{Preliminaries}
    In this section, we undertake a thorough examination of the essential geometric concepts relating to equi-affine and equi-centro-affine spaces. This pivotal understanding acts as the foundation for our subsequent analysis, enabling us to to explore the equi-affine invariants on the ellipsoid.

\subsection{Ellipsoid with equi-affine metric}
    Suppose that $N^{n+1}$ is an ellipsoid centered at the origin in $\R^{n+2}$, and under an equi-affine transformation, it may be expressed in the form given by \eqref{Ellipsoid}.
%
    A simple computation yields
    \begin{align*}
    [\mathbf{x}_r, \mathbf{x}_{\theta^1}, \cdots, \mathbf{x}_{\theta^n},\mathbf{x}] &= (-1)^{n+1} a_0\cdots a_{n+1}\sin^nr\sin^{n-1}\theta^1\cdots\sin\theta^{n-1}\\
    &=(-1)^{n+1}\sqrt{\det{\bar{g}}}.
    \end{align*}
    Consider the local embedding map
    \begin{equation*}
    \mathbf{x}:\mathbb{S}^n\mapsto M^n\hookrightarrow N^{n+1}\hookrightarrow \R^{n+2},
    \end{equation*}
    and let $D$ denote the Levi-Civita connection on $\mathbb{S}^n$.
    We typically consider $r(\theta^1,\cdots,\theta^n)$ as a function defined on the $n$-sphere $\mathbb{S}^n$.
    The derivative of $r(\theta^1,\cdots,\theta^n)$ along the direction specified by the partial derivative operator $\displaystyle \frac{\partial}{\partial\theta^i}$ is commonly denoted as $D_ir$.

    With this understanding, we can then proceed to define the local coordinate vector fields on the manifold $M^n$, utilizing these derivatives and the intrinsic structure of $\mathbb{S}^n$
    \begin{equation}\label{ek}
      e_i\triangleq \mathbf{x}_*(\frac{\partial}{\partial\theta^i})=\frac{\partial \mathbf{x}}{\partial\theta^i}=D_ir\frac{\partial}{\partial r}+\frac{\partial}{\partial\theta^i},\quad 1\leq i\leq n,
    \end{equation}
    and the outward unit normal vector of $M^n$ with respect to the metric $\bar{g}$
    \begin{equation}\label{nu}
     \nu=\frac{1}{v}\left(\frac{\partial}{\partial r}-\lambda^{-2}(r)D^jr\frac{\partial}{\partial\theta^j}\right),
    \end{equation}
    where $\lambda(r)=\sin r$, $D^jr=\sigma^{ij}D_ir$, $(\sigma^{ij})$ is the inverse matrix of the metric tensor
    \begin{equation*}
        \left(\sigma_{ij}\right)=
        \left(
          \begin{array}{cccc}
            1& 0 & \cdots & 0 \\
            0 & \sin^2\theta^1 & \cdots & 0 \\
            \cdots & \cdots & \cdots & \cdots \\
           0 & 0 & \cdots & \sin^2\theta^1\cdots\sin^2\theta^{n-1} \\
          \end{array}
        \right)
    \end{equation*}
    of $\mathbb{S}^n$.
    In addition, $v=\varrho\sqrt{1+\lambda^{-2}(r)|Dr|^2}$,  and $|\cdot|$ is the norm with respect to the metric $\sigma_{ij}$.

    It is observed that, according to \eqref{ECA}
    \begin{equation*}
    \bar{\Gamma}^k_{ij}=\hat{\Gamma}^k_{ij},\; \bar{\Gamma}^0_{ij}=-\lambda\lambda'\sigma_{ij},\;\bar{\Gamma}^k_{0i}=\frac{\lambda'}{\lambda}\delta^k_i,\;
    \bar{\Gamma}^0_{0i}=\bar{\Gamma}^k_{00}=\bar{\Gamma}^0_{00}=0,
    \end{equation*}
    where $\hat{\Gamma}^k_{ij}$ denote the Christoffel symbols of $\mathbb{S}^n$ with respect to the tangent basis $\displaystyle \{\frac{\partial}{\partial \theta^i}\}, \; i=1,\cdots,n$ and $\bar{\Gamma}^{\gamma}_{\alpha\beta}$ denote the Christoffel symbols of $N^{n+1}$ with respect to the metric $\bar{g}$. Now direct evaluation reveals
    \begin{align*}
      [e_1,e_2,\cdots, e_n,\nu,\mathbf{x}] &=\frac{1}{v}[\mathbf{x}_1,\cdots, \mathbf{x}_n, \mathbf{x}_r, \mathbf{x}]
      \left|
        \begin{array}{ccccc}
           1 & 0 & \cdots & -\frac{r^1}{\lambda^2} & 0 \\
           0 & 1 & \cdots &  -\frac{r^2}{\lambda^2}&0 \\
           \cdots & \cdots & \cdots & \cdots & \cdots \\
           r_1 & r_2 & \cdots & 1 & 0 \\
           0 & 0 & \cdots & 0 & 1 \\
         \end{array}
      \right|\\
      &= \frac{v}{\varrho^2}[\mathbf{x}_1,\cdots, \mathbf{x}_n, \mathbf{x}_r, \mathbf{x}].
    \end{align*}
    By comparison, we may find
    \begin{equation*}
      \bar{\nabla}_{e_i}e_j=\Gamma^k_{ij}e_k+h_{ij}\nu,\qquad \bar{\nabla}_{e_i}\nu=S^k_ie_k,
    \end{equation*}
   and the metric $g=g_{ij}d\theta^id\theta^j$ of $M^n$ induced from $\bar{g}$ is
    \begin{equation*}
       g_{ij}=\bar{g}(e_i,\;e_j)=\varrho^2r_ir_j+\varrho^2\lambda^2(r)\sigma_{ij}.
    \end{equation*}
    Note that in this paper subscripts after a semicolon $``{}_{,}"$ are used to denote covariant derivatives with respect to the induced equi-affine metric $g_{ij}$. Unless otherwise noted, we raise and lower indices using the metric $g_{ij}$.

    Through direct computations, we arrive at
    \begin{equation*}
       [e_1,e_2,\cdots, e_n,\nu,\mathbf{x}] = -\varrho\sqrt{\det(g_{ij})}
    \end{equation*}
    and
    \begin{align*}
       h_{ij}=g(\bar{\nabla}_{e_i}e_j,\nu)
       =\frac{\varrho^2}{v}\left(r_{i;j}-\lambda(r)\lambda'(r)\sigma_{ij}-2\frac{\lambda'(r)}{\lambda(r)}r_ir_j\right),
    \end{align*}
    where $r_{i;j}$ denotes the second covariant derivative of $r$ with respect to the Christoffel symbols $\hat{\Gamma}^k_{ij}$.

    Let $S_r$ denote the $r$-th elementary symmetric function of the the eigenvalues
    $k_1,\cdots,k_n$ of $h$, specifically:
    \begin{align*}
        S_0=1,\;S_1=k_1+\cdots+k_n,\;\cdots, \; S_n=k_1\cdots k_n.
    \end{align*}
    The Newton transformation $T_r$ are then defined inductively in the following manner
    \begin{align*}
      T_0{}^i_j=\delta^i_j,\; T_{r+1}{}^i_j=S_{r+1}\delta^i_j-T_{r}{}^{ik}h_{kj},\; r=0,1,\cdots, n-1,
    \end{align*}
    where $\delta^i_j$ is the Kronecker delta function.

\subsection{Equi-centro-affine hypersurfaces in ellipsoid}\label{sec-eca}
Let $N^{n+1}$ be an ellipsoid, centered at the origin in $\R^{n+2}$, equipped with the equi-affine metric that is precisely  specified in \eqref{Ellipsoid} and \eqref{ECA}.
Given the map
    \begin{equation*}
    \mathbf{x}:\mathbb{S}^n\mapsto M^n\hookrightarrow N^{n+1},
    \end{equation*}
    where $M^n$ is smoothly  immersed in the ellipsoid $N^{n+1}$, and the metric on $M^n$ is the equi-centro-affine metric as defined in \eqref{Metricdef}. Utilizing the relations established in \eqref{Gdef} and \eqref{Metricdef}, we proceed with the following derivation
  \begin{align*}
    G_{ij} &= [e_1,\ldots,e_n, \mathbf{x}, \bar{\nabla}_{e_i}e_j] \\
           &= h_{ij}[e_1,e_2,\ldots,e_n,\mathbf{x},\nu] = h_{ij}\varrho\sqrt{\det(g_{ij})}.
\end{align*}
Moreover,
\begin{align*}
    \tilde{g}_{ij} &= \varrho^{\frac{2}{n+2}}(\det(g_{ij}))^{\frac{1}{n+2}}(\det(h_{ij}))^{-\frac{1}{n+2}}h_{ij} = \varrho^{\frac{2}{n+2}}S_n^{-\frac{1}{n+2}}h_{ij},
\end{align*}
and
\begin{align*}
    \det(\tilde{g}_{ij}) &= \varrho^{\frac{2n}{n+2}}S_n^{\frac{2}{n+2}}\det(g_{ij}).
\end{align*}
Throughout the paper, we consistently assume that $S_n \neq 0$ on $M^n$, and further, when $n$ is an even number, we impose the condition $S_n > 0$. Consequently, the equi-centro-affine invariant area function with respect to the metric $\tilde{g}$ may be formulated as
\begin{align*}
    \Sigma_{eca} = \varrho^{\frac{n}{n+2}}\int_M S_n^{\frac{1}{n+2}}\sqrt{\det(g_{ij})}\,d\mu.
\end{align*}

\section{The variational formulas}
\subsection{The proof of Theorem \ref{thm-fv}}\label{sec-VF}
Suppose that $N^{n+1} \subset \mathbb{R}^{n+2}$ is an ellipsoid endowed with the equi-affine metric as specified in \eqref{Ellipsoid} and \eqref{ECA}.
Let $\mathbf{x}(\cdot,t): \mathbb{S}^n \rightarrow M^n \hookrightarrow N^{n+1}$
represent a smooth one-parameter family of hypersurface immersions evolving in $N^{n+1}$. We aim to consider the invariant hypersurface flow in $N^{n+1}$
\begin{equation}\label{ihf}
\frac{\partial \mathbf{x}}{\partial t} = W^k e_k + U\nu,
\end{equation}
where $W^k$ are some $(1,0)$ tensors, $U$ is equi-centro-affine invariant, and $e_k$ and $\nu$ are defined in \eqref{ek} and \eqref{nu}.
Here, the term ``invariant hypersurface flow" refers to motions of the hypersurface $\mathbf{x}$ governed by \eqref{ihf},
which possess the property of being invariant under the action of the equi-centro-affine transformation group in $\mathbb{R}^{n+2}$.
    In view of the computation in \cite{rei}, we have
    \begin{align*}
        \frac{\partial g_{ij}}{\partial t}&\;=\;-2Uh_{ij}+W_{j,i}+W_{i,j},\\
        \frac{\partial h_{ij}}{\partial t}&\;=\;U_{,ij}-Uh_{im}h^m_j+W^k_{,i}h_{kj}+W^k_{,j}h_{ki}+W^kh_{kij}+\frac{U}{\varrho^2}g_{ij},\\
        \frac{\partial \nu}{\partial t}&\;=\;-(U_i+W^lh_{il})g^{ij}e_j.
    \end{align*}
    Let $g = \det(g_{ij})$. It is easy to verify that
    \begin{align*}
      \frac{1}{g}\frac{\partial g}{\partial t}=-2\left(S_1U-W^j_{,j}\right),
    \end{align*}
    and
    \begin{eqnarray}\label{srt}
    \begin{aligned}
      \frac{\partial S_r}{\partial t}=&U(S_1S_r-(r+1)S_{r+1})+T^{ij}_{r-1}U_{,ij}\\
                    &\qquad\quad+S_{r,j}W^j+\frac{U}{\varrho^2}(n-r+1)S_r.
    \end{aligned}
    \end{eqnarray}

Consider the deformations $\mathbf{x}$ in \eqref{ihf} which leave $\partial M^n$ strongly fixed in the sense that both $U$ and its gradient vanish on $\partial M^n$.
It follows that if $M^n$ is compact and $\partial M^n$ is empty, then there is no restriction on the deformation.
The formula for the first variation with fixed boundary of the equi-centro-affine area integral for a hypersurface in an ellipsoid is
     \begin{eqnarray}\label{First-V}
     \begin{aligned}
       \frac{d}{dt}\Sigma_{eca}=&\varrho^{\frac{n}{n+2}}\int_M\frac{U}{n+2}\Bigg(T_{n-1}{}^{ij}\nabla_{e_i}\nabla_{e_j}\left(S_n^{-\frac{n+1}{n+2}}\right)\\
                           &\qquad\quad +S_n^{-\frac{n+1}{n+2}}\left(\frac{1}{\varrho^2}S_{n-1}-(n+1)S_nS_1\right)\Bigg)dV.
 \end{aligned}
 \end{eqnarray}
Hence, we have successfully completed the proof of Theorem \ref{thm-fv}.
\subsection{The second variational formula}\label{sec-vf2}
    Let us direct our concentration towards the second variation.
    Firstly, we may obtain
    \begin{align}\label{Gt}
         \frac{\partial \Gamma^k_{ij}}{\partial t}=g^{kl}(U_lh_{ij}-U_jh_{il}-U_ih_{jl}-Uh_{ijl})+W^k_{,ij}-W^mR^k_{ijm}.
    \end{align}
    In \cite{rei}, it is explicitly stated that $T_n=0$, $T_{r+1}=S_{r+1}I-AT_r$, and $T_rA=AT_r$, which implies
    \begin{eqnarray}\label{Tt}
     \begin{aligned}
      \frac{\partial T_{n-1}^{ij}}{\partial t}\;=&\;U_{,mk}\left(T_{n-1}^{ij}b^{mk}-T_{n-1}^{mi}b^{jk}\right)\\
      &\;\quad+U\left(S_1T_{n-1}^{ij}+S_ng^{ij}+\frac{1}{\varrho^2}T_{n-2}^{ij}\right)\\
      &\;\qquad +W^kT_{n-1,k}^{ij}-\left(W^j_{,k}T^{ki}_{n-1}+W^i_{,k}T^{kj}_{n-1}\right),
     \end{aligned}
    \end{eqnarray}
    where $A$ is the matrix $(h_{ij})$, and $(b^{ij})$ is the inverse matrix of $(h_{ij})$.
    Upon utilizing \eqref{Gt}, \eqref{Tt}, and \eqref{First-V}, and engaging in a non-trivial calculation process (detailed in Appendix A for completeness), we are able to formulate and present the following theorem.
    \begin{thm}\label{thm-second}
    The second variational formula at the critic point $t=t_0$ is
    \begin{eqnarray}\label{second-v}
    \begin{aligned}
    \frac{d^2}{dt^2}\Sigma_{eca}(t_0)=&\frac{\varrho^{\frac{n}{n+2}}}{n+2}\int_MU\bigg(f^{mkij}U_{,mkij}+f^{mkj}U_{,mkj}\\
      &\qquad\qquad +f^{mk}U_{,mk}+f^{m}U_{,m}+fU\bigg)dV,
     \end{aligned}
    \end{eqnarray}
    where
    \begin{align*}
       f^{mkij}=&-\frac{n+1}{n+2}S_{n}^{-\frac{2n+3}{n+2}}T^{ij}_{n-1}T^{mk}_{n-1},\\
       f^{mkj}=&-\frac{2(n+1)}{n+2}\left(S_n^{-\frac{2n+3}{n+1}}T^{mk}_{n-1}\right)_{,i}T^{ij}_{n-1},\\
       f^{mk}=&\left(S_n^{-\frac{n+1}{n+2}}\right)_{,ij}\left(T^{ij}_{n-1}b^{mk}-T^{mi}_{n-1}b^{jk}\right)-\frac{n+1}{n+2}\left(S_n^{-\frac{2n+3}{n+2}}T^{mk}_{n-1}\right)_{,ij}T^{ij}_{n-1}\\
       &\qquad -\frac{n+1}{n+2}S_n^{-\frac{2n+3}{n+2}}T^{mk}_{n-1}\left(\frac{2}{\varrho^2}S_{n-1}-nS_1S_n\right)\\
       &\qquad\qquad  +S_n^{-\frac{n+1}{n+2}}\left(\frac{1}{\varrho^2}T^{mk}_{n-1}-(n+1)(T^{mk}_0S_n+S_1T^{mk}_{n-1})\right),\\
       f^{m}=&-\frac{2(n+1)}{n+2}\left(S_n^{-\frac{2n+3}{n+2}}(S_1S_n+\frac{1}{\varrho^2}S_{n-1})\right)_{k}T^{mk}_{n-1}\\
       &\qquad\qquad+(n-2)(n+1)g^{mk}\left(S_n^{\frac{1}{n+2}}\right)_{k},\\
       f=&\left(S_n^{-\frac{n+1}{n+2}}\right)_{,ij}\left(S_1T^{ij}_{n-1}+S_ng^{ij}+\frac{1}{\varrho^2}T^{ij}_{n-2}\right)+g^{kl}h_{lij}\left(S_n^{-\frac{n+1}{n+2}}\right)_kT^{ij}_{n-1}\\
       &\quad-\frac{n+1}{n+2}\left(S_n^{-\frac{2n+3}{n+2}}(S_1S_n+\frac{1}{\varrho^2}S_{n-1})\right)_{,ij}T^{ij}_{n-1}\\
       &\qquad\;-\frac{n+1}{n+2}S_n^{-\frac{2n+3}{n+2}}\left(\frac{1}{\varrho^4}S_{n-1}^2-(n+1)S_1^2S_{n}^2-\frac{n}{\varrho^2}S_1S_{n-1}S_n\right)\\
       &\qquad\qquad+S_n^{-\frac{n+1}{n+2}}\bigg(\frac{2}{\varrho^4}S_{n-2}-\frac{n}{\varrho^2}S_n+2(n+1)S_2S_n\\
       &\qquad \qquad\quad    -\frac{n(n+1)}{\varrho^2}S_n-2(n+1)S_1^2S_n-\frac{n}{\varrho^2}S_1S_{n-1}\bigg).
    \end{align*}
    \end{thm}
\subsection{Equi-centro-affine curves on ellipsoid}
  For an ellipsoid $N^2$ in $\R^3$, an equi-affine transformation can be applied to transform it into a sphere centered precisely at the origin, with a radius denoted by $R$. Since $R$ and $\varrho$ differ only by a constant scaling factor,
  we henceforth adopt the convention of substituting $R$ for $\varrho$  in the equi-affine metric \eqref{ECA} that characterizes the ellipsoid $N^2$.

Consequently, when a curve $\mathbf{x}(p)$ is situated on this transformed spherical surface $\mathbb{S}^2$, the quantity $S_1=k_g$ directly corresponds to the geodesic curvature of the curve $\mathbf{x}$.
This substitution simplifies our analysis and allows us to directly relate the geometric properties of the curve to those of the equivalent spherical surface.

  On the spherical surface $\mathbb{S}^2$, a vector field $\mathbf{J}$  is designated as a Killing vector field along $\mathbf{x}$ if and only if it satisfies e following conditions (for a detailed discussion, refer to \cite{agm,ls}):
  \begin{eqnarray}\label{killing1}
  \begin{aligned}
    &\langle \nabla_\mathbf{T}\mathbf{J},\mathbf{T}\rangle=0,\\
    &\langle \nabla^2_\mathbf{T}\mathbf{J},\boldsymbol{\epsilon}\rangle+\frac{1}{R^2}\langle \mathbf{J},\boldsymbol{\epsilon}\rangle=0,
  \end{aligned}
  \end{eqnarray}
  where $\mathbf{T}=\left(e/{\sqrt{g(e,e)}}\right)$ represents the unit tangent vector to the curve, with $e$ is defined as in \eqref{ek}, and $\boldsymbol{\epsilon}=\nu$ in \eqref{nu}.
  \subsubsection{The variational formulations of the equi-centro-affine arc length}
  By utilizing \eqref{First-V}, the first equi-centro-affine variational formula for the curves on the sphere $\mathbb{S}^2$ may be elegantly expressed as
  \begin{align}\label{1-vf}
     {}^{ec}\hspace{-1mm}L'(t)=\frac{1}{3}R^{1/3}\int^{s_2}_{s_1}U\left(B_{ss}-2B^{-2}+\frac{B}{R^2}\right)ds,
  \end{align}
  where $B=(\kappa_g)^{-2/3}$, and $ds=\sqrt{g(e,e)}dp$.
  At the critic point $t_0$, $B_{ss}-2B^{-2}+BR^{-2}=0$ holds. Specifically, according to \eqref{second-v}, the second variational formula is given by
  \begin{align*}
     {}^{ec}\hspace{-1mm}L''(t_0)=\frac{1}{3}R^{1/3}\int^{s_2}_{s_1}U\left(f_4U_{ssss}+f_3U_{sss}+f_2U_{ss}+f_1U_{s}+f_0U\right)ds,
  \end{align*}
  where
  \begin{align*}
   f_0&=\;\frac{1}{R^4}B^{\frac{5}{2}}-\frac{5}{2R^2}B^{\frac{1}{2}}B_s^2-\frac{9}{R^2}B^{-\frac{1}{2}}-2B^{-\frac{5}{2}}B_s^2+2B^{-\frac{7}{2}},\\
   f_1&=\;-\frac{10}{3R^2}B^{\frac{3}{2}}B_s+\frac{5}{3}B^{-\frac{3}{2}}B_s,\;\; f_3=\;-\frac{10}{3}B^{\frac{3}{2}}B_s, \\
  f_2&=\;-\frac{5}{2}B^{\frac{1}{2}}B_s^2+\frac{1}{3R^2}B^{\frac{5}{2}}-\frac{20}{3}B^{-\frac{1}{2}}, \;\; f_4=\;-\frac{2}{3}B^{\frac{5}{2}}.
  \end{align*}
  Integration by parts gives
  \begin{align*}
     \int^{s_2}_{s_1}f_4UU_{ssss}ds&=\;\int^{s_2}_{s_1}\left(f_4U_{ss}^2-2(f_4)_{ss}U_s^2+\frac{1}{2}(f_4)_{ssss}U^2\right)ds,\\
     \int^{s_2}_{s_1}f_3UU_{sss}ds&=\;\int^{s_2}_{s_1}\left(\frac{3}{2}(f_3)_sU_{s}^2-\frac{1}{2}(f_3)_{sss}U^2\right)ds,\\
     \int^{s_2}_{s_1}f_2UU_{ss}ds&=\;\int^{s_2}_{s_1}\left(-f_2U_{s}^2+\frac{1}{2}(f_2)_{ss}U^2\right)ds,\\
     \int^{s_2}_{s_1}f_1UU_{s}ds&=\;\int^{s_2}_{s_1}-\frac{1}{2}(f_1)_{s}U^2ds.\\
  \end{align*}
  Then one obtains
  \begin{align}\label{1-vf2}
     {}^{ec}\hspace{-1mm}L''(t_0)=\frac{1}{3}R^{1/3}\int^{s_2}_{s_1}\left(P_2U_{ss}^2+P_1U_s^2+P_0U^2\right)ds,
  \end{align}
  where
  \begin{align*}
  P_0&=\;\frac{1}{2}(f_4)_{ssss}-\frac{1}{2}(f_{3})_{sss}+\frac{1}{2}(f_2)_{ss}-\frac{1}{2}(f_1)_s+f_0,\\
   P_1&=\;-2(f_4)_{ss}+\frac{3}{2}(f_3)_s-f_2,\;\; P_2=\;f_4.
  \end{align*}
Hence, we have
  \begin{align*}
    P_0&=\;\frac{1}{R^4}B^{\frac{5}{2}}-\frac{5}{2R^2}B^{\frac{1}{2}}B_s^2-\frac{9}{R^2}B^{-\frac{1}{2}}-2B^{-\frac{5}{2}}B_s^2+2B^{-\frac{7}{2}},\\
    P_1&=\;\frac{4}{3R^2}B^{\frac{5}{2}}+\frac{10}{3}B^{-\frac{1}{2}},\;\; P_2=\;-\frac{2}{3}B^{\frac{5}{2}}.
  \end{align*}

  Solitons of the geometric flow represent self-similar solutions that maintain their shape under the evolution induced by the flow. To formulate the definition of a soliton, let $N^n$ be a $n$-dimensional Riemannian manifold with metric $g$, equipped with a Killing vector field $\mathbf{J}$ related to an isometry group $\varphi:N\times\R\to N$.
   The isometry group $\varphi$ characterizes transformations that preserve the metric $g$ and hence the geometric properties of the manifold. The relationship between $\mathbf{J}$ and $\varphi$ is given by the following differential equation and initial condition
  \begin{align*}
      \frac{d\varphi(x,t)}{dt}&=\;\mathbf{J}(\varphi(x,t)),\\
      \varphi(x,0)&=x.
  \end{align*}
  Here, $\frac{d\varphi(x,t)}{dt}$ denotes the time derivative of the point $\varphi(x,t)$ in the direction of the flow induced by the Killing vector field $\mathbf{J}$.
  The initial condition $\varphi(x,0)=x$ signifies that at time $t=0$, the isometry group leaves each point $x$ in $N^n$ unchanged.

  A curve $\mathbf{x}$ on $N^n$ is a soliton of the geometric flow \eqref{ihf} if, under the action of the isometry group $\varphi$ parameterized by the flow induced by $\mathbf{J}$, the curve evolves in such a way that its shape remains constant up to reparametrization. This property ensures that the soliton maintains its geometric characteristics throughout the evolution process.
  \begin{thm}\label{thm-soliton}
  If $\displaystyle U=a(t)\cos\frac{2\pi s}{{}^e\hspace{-0.9mm}L(t)}+b(t)\sin\frac{2\pi s}{{}^e\hspace{-0.9mm}L(t)}$ and $W_s=\kappa_gU$ in \eqref{ihf}, then  the closed curve $\mathbf{x}$ is the soliton of the flow \eqref{ihf} if and only if  its geodesic curvature remains constant. Here ${}^e\hspace{-0.9mm}L(t)$ represents the Euclidean arc length of $\mathbf{x}$.
  \end{thm}
  \begin{proof}
     If $\kappa_g$ of the closed curve $\mathbf{x}$ is constant, we can confirm that the vector field $U\boldsymbol{\epsilon}+W\mathbf{T}$ constitutes a Killing vector field for $\mathbf{x}$ by satisfying the conditions outlined in \eqref{killing1}.
  \end{proof}
  \begin{thm}
    A closed curve with a geodesic curvature of $\kappa_g=\frac{\sqrt{2}}{2R}$ is an equi-centro-affine maximal curve that lies on a sphere, where $R$ represents the radius of the sphere.
  \end{thm}
  \begin{proof}
    Given that the geodesic curvature $\kappa_g=\frac{\sqrt{2}}{2R}$, it follows that the curve is a closed planar circle with a radius of $r=\frac{\sqrt{6}}{3}R$.
    Then by \eqref{1-vf} and \eqref{1-vf2}, at the critic point  $t=t_0$,  we have $${}^{ec}\hspace{-1mm}L'(t_0)=0$$
    and
    \begin{align*}
      {}^{ec}\hspace{-1mm}L''(t_0)&=\;\frac{2^{\frac{5}{6}}}{3}R^2\int^{\frac{2\sqrt{6}\pi}{3}R}_0\left(-\frac{2}{3}U_{ss}^2+3R^{-2}U_s^2-3R^{-4}U^2\right)ds\\
                     &=\;-\frac{2^{\frac{5}{6}}\sqrt{6}}{6R}\int^{2\pi}_0\left(U_{xx}^2-3U_x^2+2U^2\right)dx.
    \end{align*}
    Let
    \begin{align*}
      U=\frac{a_0}{2}+\sum^{\infty}_{m=1}(a_m\cos(mx)+b_m\sin(mx)).
    \end{align*}
    According to Theorem \ref{thm-soliton}, if the function $U$ is expressed as  $U=a\cos(x)+b\sin(x)$,  then the resulting closed curves on the sphere with a constant geodesic curvature  $\kappa_g$
    are identified as solitons of the given equation \eqref{ihf}.
    Consequently, we proceed under the assumption that
    $a_0, a_2, b_2, a_3, b_3, \cdots $ are not all zeros.
    Then
      \begin{align*}
      \int_0^{2\pi}U^2dx&=\;\frac{a_0^2}{2}\pi+\sum^{\infty}_{m=0}(a_m^2+b_m^2)\pi,\\
      \int_0^{2\pi}U^2_xdx&=\;\sum^{\infty}_{m=1}m^2(a_m^2+b_m^2)\pi,\\
      \int_0^{2\pi}U^2_{xx}dx&=\;\sum^{\infty}_{m=1}m^4(a_m^2+b_m^2)\pi,
    \end{align*}
    which imply
    \begin{align*}
      {}^{ec}\hspace{-1mm}L''(t_0)&=\;-\frac{2^{\frac{5}{6}}\sqrt{6}}{6R}\left(a_0^2\pi+\pi\sum^{\infty}_{m=1}(m^2-2)(m^2-1)(a_m^2+b_m^2)\right).
    \end{align*}
    Since $m$ is an integer,  $(m^2-2)(m^2-1)\geq 0$. Thus
    \begin{align*}
    {}^{ec}\hspace{-1mm}L''(t_0)<0.
    \end{align*}
    Thus the theorem is proved.
  \end{proof}
  \subsubsection{The area-fixed variational formulation}\label{sec-fix}
  On the sphere, let us consider a closed curve denoted by $\mathbf{x}$, which possesses an Euclidean arc length ${{}^{e}\hspace{-1mm}L}$. Furthermore, let $A$ represent the area enclosed by this curve $\mathbf{x}$. According to the Gauss-Bonnet Theorem,  we have the following fundamental relationship
  \begin{align*}
      \frac{1}{R^2}A+\int^{{}^{e}\hspace{-1mm}L}_0\kappa_gds=2\pi.
  \end{align*}
  Assuming that the motion governed by \eqref{ihf} is area-preserving, a crucial consequence arises from the conservation of enclosed area. Specifically,
  \begin{align*}
      \frac{d}{dt}\int^{{}^{e}\hspace{-1mm}L}_0\kappa_gds=0.
  \end{align*}
  According to \eqref{srt}, the above equation implies
  \begin{align*}
     \int^{{}^{e}\hspace{-1mm}L}_0 Uds=0.
  \end{align*}
  By \eqref{1-vf}, at the critic point  $t=t_0$,
  \begin{align*}
    \int^{{}^{e}\hspace{-1mm}L}_{0}U\left(B_{ss}-2B^{-2}+\frac{B}{R^2}\right)ds=0.
  \end{align*}
  Let us denote
  \begin{align*}
     C_1=\int^{{}^{e}\hspace{-1mm}L}_0 \left(B_{ss}-2B^{-2}+\frac{B}{R^2}\right) ds=\int^{{}^{e}\hspace{-1mm}L}_0 \left(\frac{B}{R^2}-2B^{-2}\right) ds.
  \end{align*}
  Clearly,
  \begin{align*}
     \int^{{}^{e}\hspace{-1mm}L}_0C_1Uds=0,
  \end{align*}
  and
  \begin{align}\label{Uint}
    \int^{{}^{e}\hspace{-1mm}L}_{0}U\left(B_{ss}-2B^{-2}+\frac{B}{R^2}-\frac{C_1}{{}^{e}\hspace{-1mm}L}\right)ds=0.
  \end{align}
  Of course,
  \begin{align*}
    \int^{{}^{e}\hspace{-1mm}L}_{0}\left(B_{ss}-2B^{-2}+\frac{B}{R^2}-\frac{C_1}{{}^{e}\hspace{-1mm}L}\right)ds=0.
  \end{align*}
  If $U=B_{ss}-2B^{-2}+\frac{B}{R^2}-\frac{C_1}{{}^{e}\hspace{-1mm}L}$, substituting it into \eqref{Uint} yields
  \begin{align*}
    B_{ss}-2B^{-2}+\frac{B}{R^2}=\frac{C_1}{{}^{e}\hspace{-1mm}L},
  \end{align*}
  which signifies that the expression $B_{ss}-2B^{-2}+\frac{B}{R^2}$ is constant.

  Let us now delve into the analysis of the stability of the variational formulation, under the assumption that $B$ remains constant.
  Under this condition, the configuration of $\mathbf{x}$ adopts the shape of a planar circle, with a radius $r$ given by $r=\frac{1}{\sqrt{R^{-2}+B^{-3}}}$.
 Furthermore, the circumference ${{}^{e}\hspace{-1mm}L}$ of this circle is simply ${{}^{e}\hspace{-1mm}L}=2r\pi$.
  \begin{thm}
     If $B$ is constant and within the range $\frac{7}{5}R^2\leq B^3 \leq 2R^2$, then $\mathbf{x}$ is an equi-centro-affine maximal curve on the sphere under under the constraint of area-preserving motions.
  \end{thm}
  \begin{proof} If $B$ is constant, then
  \begin{align*}
    {}^{ec}\hspace{-1mm}L''(t_0)&=\;\frac{1}{3}R^{\frac{1}{3}}\int^{2r\pi}_0\Big(-\frac{2}{3}B^{\frac{5}{2}}U_{ss}^2+\left(\frac{4}{3R^2}B^{\frac{5}{2}}+\frac{10}{3}B^{-\frac{1}{2}}\right)U_s^2\\
    &\qquad \qquad+\left(\frac{1}{R^4}B^{\frac{5}{2}}-\frac{9}{R^2}B^{-\frac{1}{2}}+2B^{-\frac{7}{2}}\right)U^2\Big)ds\\
    &=\;\frac{r}{3}R^{\frac{1}{3}}\int^{2\pi}_0\Big(-\frac{2}{3r^4}B^{\frac{5}{2}}U_{xx}^2+\left(\frac{4}{3r^2R^2}B^{\frac{5}{2}}+\frac{10}{3r^2}B^{-\frac{1}{2}}\right)U_x^2\\
    &\qquad\qquad+\left(\frac{1}{R^4}B^{\frac{5}{2}}-\frac{9}{R^2}B^{-\frac{1}{2}}+2B^{-\frac{7}{2}}\right)U^2\Big)dx\\
    &=\;-\frac{r}{9}R^{\frac{1}{3}}\int^{2\pi}_0\Big(\frac{1}{R^4}B^{\frac{5}{2}}\left(2U_{xx}^2-4U_x^2-3U^2\right)\\
    &\qquad\qquad +\frac{1}{R^2}B^{-\frac{1}{2}}\left(4U_{xx}^2-14U_x^2+27U^2\right)\\
    &\qquad\qquad\qquad+B^{-\frac{7}{2}}\left(2U_{xx}^2-10U_x^2-6U^2\right)\Big)dx.
  \end{align*}
  Since $\displaystyle \int^{2r\pi}_0Uds=0$, we denote $U$ by
  \begin{align*}
      U=\sum^{\infty}_{m=1}(a_m\cos(mx)+b_m\sin(mx)),
    \end{align*}
  where $a_2, b_2, a_3, b_3, \cdots $ are not all zeros.

  Thus,
  \begin{align*}
    {}^{ec}\hspace{-1mm}L''(t_0)&=\;-\frac{r}{9}R^{\frac{1}{3}}B^{-\frac{7}{2}}\sum^\infty_{m=1}\Big((2m^4-4m^2-3)\frac{B^{6}}{R^4}\\
                   &\qquad\qquad+(4m^4-14m^2+27)\frac{B^{3}}{R^2}\\
                   &\qquad\qquad\qquad+(2m^4-10m^2-6)\Big)(a_m^2+b_m^2)\pi.
  \end{align*}
  when $m>2$, we may verify $2m^4-4m^2-3>0$, $4m^4-14m^2+27>0$ and $2m^4-10m^2-6>0$.
  Conversely, if $m=1$ (or $m=2$), by $\frac{7}{5}R^2\leq B^3 \leq 2R^2$, we see
  \begin{eqnarray*}
  \begin{aligned}
  (2m^4-4m^2-3)\frac{B^{6}}{R^4}+(4m^4-14m^2+27)\frac{B^{3}}{R^2}+(2m^4-10m^2-6)\geq0\; (>0).
  \end{aligned}
  \end{eqnarray*}
  Therefore, we obtain
  \begin{align*}
    {}^{ec}\hspace{-1mm}L''(t_0)<0.
  \end{align*}
  This proves the theorem.
  \end{proof}

\section{Equi-centro-affine extremal hypersurfaces in unit sphere}\label{sec-eh}
 Li \cite{lhz} derived the isoparametric Willmore hypersurfaces. In a similar fashion, we are tasked with classifying the isoparametric equi-centro-affine extremal hypersurfaces. The subsequent lemma, which serves as a crucial foundation for our classification, was initially presented in \cite{lhz} (refer to \cite{abr,lhz, mun,sto} for further details).
\begin{lem}[\cite{lhz}]\label{lem-iso}
 Let $M$ be an $n$-dimensional compact isoparametric hypersurface (i.e. hypersurface with constant principal curvatures) in $\mathbb{S}^{n+1}(1)$. Let $k_1 > k_2 >\cdots > k_g$ be the distinct principal curvatures with multiplicities $m_1, m_2, \cdots , m_g$ ( so that $n = m_1 +m_2 + \cdots + m_g$). Then
\begin{itemize}
  \item[(1)] $g$ is either $1,2,3,4,$ or $6$.
  \item[(2)] If $g = 1$,\; $M$ is totally umbilic.
  \item[(3)] If $g = 2$,\; $\displaystyle M = \mathbb{S}^m\left(r_1\right) \times \mathbb{S}^{n-m}\left(r_2\right), \quad r_1^2+r_2^2 = 1$.
  \item[(4)] If $g = 3$,\;  $m_1 =m_2=m_3 = 2^k, \quad (k = 0,1,2,3)$.
  \item[(5)] If $g = 4$, $m_1 = m_3$ and $m_2 = m_4$. Moreover, $(m_1,m_2)=(2,2)$ or $(4,5)$, or
$m_1 + m_2 + 1$ is a multiple of $2^{\phi(m_1-1)}$. Here $\phi(l)$ is the number of integers $s$ with
$1\leq s \leq l$ and $s = 0,1,2,4 \mod 8$.
  \item[(6)] If $g = 6$, $m_1 = m_2 = \cdots = m_6 = 1$ or $2$.
\item[(7)] There exists an angle $\theta$, $\displaystyle 0 < \theta <\frac{\pi}{g}$ such that
\begin{align}\label{kalpha}
k_{\alpha}=\cot\left(\theta+\frac{\alpha-1}{g}\pi\right),\quad \alpha=1,2,\cdots g.
\end{align}
\end{itemize}
\end{lem}
Initially, leveraging \eqref{extremal}, we deduce the ensuing result.
\begin{lem}
 Let $M$ be an $n$-dimensional equi-centro-affine extremal hypersurface in $\mathbb{S}^{n+1}$ characterized by a constant value of $S_n$. Then, the following fundamental identity holds:
 \begin{align}\label{cons-sn}
  S_{n-1} - (n + 1)S_nS_1 = 0.
 \end{align}
\end{lem}
\begin{thm}\label{thm-iso-eca}
 Let $M$ be an $n$-dimensional compact isoparametric equi-centro-affine extremal hypersurface in $\mathbb{S}^{n+1}$. Then the classification of such hypersurfaces, based on their number $g$ of distinct principal curvatures, is as follows:
 \begin{itemize}
  \item[(i)] If $g = 1$, $M$ is totally umbilic. Specifically, the principal curvatures satisfy $k^2 = \frac{1}{n+1}$ and the radius is $\sqrt{\frac{n+1}{n+2}}$.
  \item[(ii)] If $g = 2$, $M$ is given by
   \begin{align}
     M = \mathbb{S}^m\left(\sqrt{\frac{m+1}{n+2}}\right) \times \mathbb{S}^{n-m}\left(\sqrt{\frac{n+1-m}{n+2}}\right),
   \end{align}
   where $1 \leq m \leq n-1$. In particular, if $n = 2m$, and $m$ should be an even number.
  \item[(iii)] If $g = 3$, the principal curvatures $k_1, k_2, k_3$ satisfy:
   \begin{align*}
    k_1^2k_2^2k_3^2 &= \frac{1}{n+1}, \\
    k_1^2 + k_2^2 + k_3^2 &= \frac{3(2n+5)}{n+1}, \\
    k_1^2k_2^2 + k_2^2k_3^2 + k_1^2k_3^2 &= \frac{3(3n+5)}{n+1},
   \end{align*}
   and the dimension $n$ can only take the values $\{3, 6, 12, 24\}$.
  \item[(iv)] If $g = 4$, the principal curvatures $k_1, k_2, k_3, k_4$ satisfy:
   \begin{align*}
    k_1 &= 1 + \sqrt{2}, & k_2 &= \sqrt{2} - 1, \\
    k_3 &= 1 - \sqrt{2}, & k_4 &= -(1 + \sqrt{2}),
   \end{align*}
   and the dimension $n$ is uniquely determined to be $8$.
  \item[(v)] If $g = 6$, the principal curvatures $k_1, \ldots, k_6$ satisfy:
   \begin{align*}
    k_1 &= 2 + \sqrt{3}, & k_2 &= 1, & k_3 &= 2 - \sqrt{3}, \\
    k_4 &= -(2 - \sqrt{3}), & k_5 &= -1, & k_6 &= -(2 + \sqrt{3}),
   \end{align*}
   and the dimension $n$ is fixed at $12$.
 \end{itemize}
\end{thm}
\begin{proof}
(i) If $g=1$, we have
\begin{align*}
   nk_1^{n-1}-n(n+1)k_1^{n+1}=0,
\end{align*}
which implies
\begin{align*}
  k_1^2=\frac{1}{n+1}.
\end{align*}

(ii) If $g = 2$, let distinct principal curvatures are $k_1$ (multiplicity $m$) and $k_2$ (multiplicity $n-m$).
Then by (3) of Lemma \ref{lem-iso}, \eqref{kalpha} and \eqref{cons-sn}, we have
\begin{gather*}
  1+k_1k_2=0,\\
  (n-m)k_1^mk_2^{n-m-1}+mk_1^{m-1}k_2^{n-m}-(n+1)k_1^mk_2^{n-m}\left(mk_1+(n-m)k_2\right)=0.
\end{gather*}
It follows that
\begin{align*}
  k_1^2=\frac{n+1-m}{m+1}.
\end{align*}
\begin{align*}
 M&\;=\;\mathbb{S}^m\left(\frac{1}{\sqrt{1+k_1^2}}\right)\times\mathbb{S}^{n-m}\left(\frac{1}{\sqrt{1+\frac{1}{k_1^2}}}\right)\\
 &\;=\;\mathbb{S}^m\left(\sqrt{\frac{m+1}{n+2}}\right)\times\mathbb{S}^{n-m}\left(\sqrt{\frac{n+1-m}{n+2}}\right),
\end{align*}
where $\displaystyle 1\leq m\leq n-1$.

(iii) If $g=3$, by (4) of Lemma \ref{lem-iso}, $m_1=m_2=m_3:=m$, $n=3m$.
From (7) of Lemma \ref{lem-iso}, we may deduce
\begin{align*}
 k_1=\cot\theta, \quad k_2=\frac{k_1-\sqrt{3}}{1+\sqrt{3}k_1},\quad k_3=\frac{k_1+\sqrt{3}}{1-\sqrt{3}k_1}.
\end{align*}
Then
\begin{align}\label{k1eqg3}
  (3m+1)k_1^6-3(6m+5)k_1^4+3(9m+5)k_1^2-1=0.
\end{align}
It is direct to confirm that $k_2$ and $k_3$ are likewise solutions of \eqref{k1eqg3}. Hence,
\begin{align*}
k_1^2k_2^2k_3^2&\;=\;\frac{1}{n+1},\\ k_1^2+k_2^2+k_3^2&\;=\;\frac{3(2n+5)}{n+1},\\ k_1^2k_2^2+k_2^2k_3^2+k_1^2k_3^2&\;=\;\frac{3(3n+5)}{n+1},
\end{align*}
where $n\in\{3,6,12,24\}$.

(iv) If $g=4$, according to (4) and (7) of Lemma \ref{lem-iso}, we obtain $m_1=m_3,\;m_2=m_4$.
Furthermore, the specific relationships between the principal curvatures can be expressed as
\begin{align*}
k_2=\frac{k_1-1}{k_1+1},\quad k_3=-\frac{1}{k_1}, \quad k_4=-\frac{k_1+1}{k_1-1}.
\end{align*}
Denote
\begin{align*}
 A=k_1+k_3,\quad B=k_2+k_4.
\end{align*}
Then
\begin{align*}
S_1&\;=\;m_1A+m_2B;\\
S_n&\;=\;k_1^{m_1}k_2^{m_2}k_3^{m_3}k_4^{m_4}=(-1)^{m_1+m_2},\\
S_{n-1}&\;=\;S_n\left(\frac{m_1}{k_1}+\frac{m_2}{k_2}+\frac{m_3}{k_3}+\frac{m_4}{k_4}\right)=(-1)^{m_1+m_2+1}S_1,
\end{align*}
and
\begin{align*}
S_{n-1}-(n+1)S_1S_n=S_1\left((-1)^{m_1+m_2+1}-(n+1)(-1)^{m_1+m_2}\right)=0.
\end{align*}
Then we attain
\begin{align*}
   S_1=m_1A+m_2B=0.
\end{align*}
It is noteworthy that  $AB=-4$, and from this, we derive
\begin{align*}
  A^2=\frac{4m_2}{m_1}.
\end{align*}
The fact that $n=2(m_1+m_2)$ is an even number necessitates the condition  $S_n=(-1)^{m_1+m_2}>0$, which in turn implies that $m_1+m_2$ must an even number.
Applying (4) of Lemma \ref{lem-iso}, we deduce that $m_1=m_2=2$, resulting in  $n=8$.
Subsequently, the specific values of the principal curvatures can be determined as
\begin{align*}
    k_1=1+\sqrt{2}, \; k_2=\sqrt{2}-1,\;k_3=1-\sqrt{2},\;k_4=-(1+\sqrt{2}),
  \end{align*}

(v) If $g=6$,  by (6) of Lemma \ref{lem-iso}, we deduce that  $m_1=\cdots=m_6$ can only take values of $1$ or $2$. Further, based on part (7) of the same lemma, we have the following expressions
\begin{gather*}
k_1=\cot\theta, \quad k_2=\frac{\sqrt{3}k_1-1}{k_1+\sqrt{3}}, \quad k_3=\frac{k_1-\sqrt{3}}{1+\sqrt{3}k_1},\\
k_4=-\frac{1}{k_1},\quad k_5=-\frac{1}{k_2} \quad k_6=-\frac{1}{k_3}.
\end{gather*}
Using these values, we can express the sums as follows
\begin{align*}
S_1=m_1\sum_{i=1}^6k_i,\quad S_n=1, \quad S_{n-1}=-m_1S_1,
\end{align*}
and
\begin{align*}
  S_{n-1}-(n+1)S_1S_n=S_1(-m_1-(6m_1+1))=0.
\end{align*}
Then $S_1=0$ results in
\begin{gather*}
   k_1=2+\sqrt{3},\;k_2=1,\;k_3=2-\sqrt{3},\\ k_4=-(2-\sqrt{3}),\; k_5=-1,\;k_6=-(2+\sqrt{3}).
\end{gather*}
We exclude the case where  $n=6$ due to the fact that $S_6=-1<0$.

Accordingly, the proof is now comprehensively established, thereby concluding our argument.
\end{proof}
\begin{rem} Among these isoparametric equi-centro-affine extremal hypersurfaces outlined in Theorem \ref{thm-iso-eca}, the hypersurfaces specified in (ii) with the condition $2m=n$, as well as those in (iv) and (v)
are Euclidean isoparametric minimal hypersurfaces in $\mathbb{S}^{n+1}(1)$.
In addition, the hypersurfaces in (iv) and (v)
are isoparametric Willmore hypersurfaces.
\end{rem}
\begin{thm}
Let $M^n$ be a compact equi-centro-affine extremal hypersurface
in  unit sphere $\mathbb{S}^{n+1}$. Consequently, the ensuing conclusions hold true.
\begin{itemize}
  \item If $\displaystyle 0<k_i\leq\frac{1}{\sqrt{n+1}},\; i=1\;,\cdots,\;n$, then $M^n$ is a totally umbilic hypersphere with $\displaystyle k_i=\frac{1}{\sqrt{n+1}},\;  i=1,\cdots,n$.
  \item If $\displaystyle k_i\geq\frac{1}{\sqrt{n+1}},\; i=1\;,\cdots,\;n$, then $M^n$ is a totally umbilic hypersphere with $\displaystyle k_i=\frac{1}{\sqrt{n+1}},\;  i=1,\cdots,n$.
\end{itemize}
\end{thm}
\begin{proof}
 Pursuant to Proposition B as outlined in \cite{rei}, we see the Newton tensors are divergence-free, explicitly denoted as $\displaystyle T_r{}^{ij}_{,j}=0$.
 Thus,
 \begin{align*}
  0&\;=\;\int_M T_{n-1}{}^{ij}\nabla_{e_i}\nabla_{e_j}\left(S_n^{-\frac{n+1}{n+2}}\right)+S_n^{-\frac{n+1}{n+2}}\left(S_{n-1}-(n+1)S_nS_1\right) d\mu_{M}\\
  &\;=\;\int_M S_n^{-\frac{n+1}{n+2}}\left(S_{n-1}-(n+1)S_nS_1\right) d\mu_{M}\\
  &\;=\;\int_M S_n^{\frac{1}{n+2}}\sum_{i=1}^n\frac{1-(n+1)k_i^2}{k_i} d\mu_{M}.
 \end{align*}
 Since $S_n\neq0$ throughout $M^n$, the aforementioned formula yields the anticipated outcomes as desired.
\end{proof}
\section{The closed equi-centro-affine extremal curves on sphere}
In this section, we undertake a comprehensive classification of closed curves residing on the sphere that fulfill the intricate equation
\begin{align}
    B_{ss}-2B^{-2}+R^{-2}B=\frac{C_1}{2}, \label{BC-eq1}
  \end{align}
  where $\displaystyle B=\kappa_g^{-\frac{2}{3}}$ and $C_1$ is constant.  It is noteworthy that this equation holds a pivotal position in characterizing the critical points of equi-centro-affine arc length, subject to the constraint of a fixed enclosed area on the spherical surface. By categorizing these curves, we gain an understanding of their geometric properties, their relationship to arc length metrics, and how they are influenced by the inherent constraints of the spherical geometry.
  \begin{defn} On the sphere $\mathbb{S}^2$, a curve  whose geodesic curvature satisfies \eqref{BC-eq1} is called  {\it generalized equi-centro-affine extremal} curve.
  In particular, if $C_1\equiv0$, it is called {\it equi-centro-affine extremal} curve.
  \end{defn}
  Upon integrating \eqref{BC-eq1}, we arrive at the subsequent equation, which includes the integration constant  $C_2$.
  \begin{align}
      B_s^2=C_2-\frac{B^2}{R^2}-4B^{-1}+C_1B. \label{BC-eq2}
  \end{align}
  Notice that, from the preceding equation, we derive the following crucial inequality
  \begin{align*}
      C_2\geq\left(\frac{B}{R}-\frac{RC_1}{2}\right)^2+4B^{-1}-\frac{R^2C_1^2}{4}>-\frac{R^2C_1^2}{4}.
  \end{align*}
Given the assumption that $B$ is not constant, we now proceed to analyze the solutions of the polynomial equation $B^3-R^2C_1B^2-R^2C_2B+4R^2=0$.
  \begin{prop}\label{prop-bd}
    The polynomial equation $B^3-\mu B^2-\lambda B+4R^2=0$, where $\mu=R^2C_1$ and $\lambda=R^2C_2$,
    admits two distinct positive solutions provided that either of the following two conditions is fulfilled:
   \begin{itemize}
     \item[(a)] If $\displaystyle \mu\geq-3(2R)^{\frac{2}{3}}$, then
     \begin{align*}
      \lambda>\frac{-\mu^2-(Y_1^{1/3}+Y_2^{1/3})}{12},
     \end{align*}
     where
      $\displaystyle
         Y_{1,2}=\mu^6+2160R^2\mu^3-93312R^4\pm48\sqrt{3}R\sqrt{(\mu^3+108R^2)^3}.
      $
     \item[(b)] If $\displaystyle \mu\leq-3(2R)^{\frac{2}{3}}$, then
     \begin{align*}
     \lambda>\frac{-\mu^2}{12}+\frac{\sqrt{\mu(\mu^3-864R^2)}}{6}\cos\frac{\vartheta-\pi}{3},
   \end{align*}
     where $\displaystyle \vartheta=\arccos\frac{-\mu^6-2160R^2\mu^3+93312R^4}{\mu(864R^2-\mu^3)\sqrt{\mu(\mu^3-864R^2)}}$.
   \end{itemize}
  \end{prop}
  \begin{proof}
   Given that $B^3-\mu B^2-\lambda B+4R^2=0$ possesses two distinct positive solutions, upon differentiating and analyzing the resulting quadratic equation $3B^2-2\mu B-\lambda=0$,
   it becomes evident that this quadratic cannot admit two negative roots. Consequently,
   if $\mu\leq0$,  then it is imperative that $\lambda>0$ holds true.
   Concurrently,  the fact that $B^3-\mu B^2-\lambda B+4R^2=0$ admits three distinct real solutions directly implies the discriminant
   \begin{align*}
     \Delta=-3(4\lambda^3+\mu^2\lambda^2+72\mu R^2\lambda+16\mu^3R^2-432R^4)<0.
   \end{align*}

   Now, let us shift our focus to considering the inequality  $4\lambda^3+\mu^2\lambda^2+72\mu R^2\lambda+16\mu^3R^2-432R^4>0$.
   To gain insight, we commence by analyzing the roots of $4\lambda^3+\mu^2\lambda^2+72\mu R^2\lambda+16\mu^3R^2-432R^4=0$. Direct computations yield the coefficients
   \begin{align*}
     \bar{A}&=\;\mu(\mu^3-864R^2),\\
     \bar{B}&=\;-72R^2(7\mu^3-216R^2), \\
     \bar{C}&=\;-48R^2\mu^2(\mu^3-135R^2)
   \end{align*}
   along with the discriminant
   \begin{align*}
     \Delta_{\Delta}=\bar{B}^2-4\bar{A}\bar{C}=192R^2(\mu^3+108R^2)^3.
   \end{align*}
   It is immediately apparent that $\bar{A}$ and $\bar{B}$ cannot simultaneously be zero. Consequently, we proceed with the following three distinct cases based on the behavior of these coefficients and the discriminant.

   {\bf Case 1.} If $\Delta_{\Delta}<0$, which translates to $\mu^3<-108R^2$, then the equation $4\lambda^3+\mu^2\lambda^2+72\mu R^2\lambda+16\mu^3R^2-432R^4=0$ possesses three distinct real roots, denoted as $\lambda_1, \lambda_2$ and $\lambda_3$,
   satisfying the ordering
   $\lambda_1<\lambda_2<0<\lambda_3$. By previous statement (if $\mu<0$, then $\lambda>0$),
   $4\lambda^3+\mu^2\lambda^2+72\mu R^2\lambda+16\mu^3R^2-432R^4>0$ yields
   \begin{align*}
     \lambda>\lambda_3=\frac{-\mu^2}{12}+\frac{\sqrt{\mu(\mu^3-864R^2)}}{6}\cos\frac{\theta-\pi}{3},
   \end{align*}
   where $\displaystyle \theta=\arccos\frac{-\mu^6-2160R^2\mu^3+93312R^4}{\mu(864R^2-\mu^3)\sqrt{\mu(\mu^3-864R^2)}}$.

   {\bf Case 2.} If $\Delta_{\Delta}<0$, which signifies that $\mu^3>-108R^2$, then the equation $4\lambda^3+\mu^2\lambda^2+72\mu R^2\lambda+16\mu^3R^2-432R^4=0$ admits a unique real root, labeled $\lambda_1$ , accompanied by a conjugate pair of complex roots,  $\lambda_2$ and $\lambda_3$.
   The real root can be expressed as
   \begin{align*}
     \lambda_1=\frac{-\mu^2-(Y_1^{1/3}+Y_2^{1/3})}{12},
   \end{align*}
   where
   \begin{align*}
     Y_{1,2}=\mu^6+2160R^2\mu^3-93312R^4\pm48\sqrt{3}R\sqrt{(\mu^3+108R^2)^3}.
   \end{align*}
   Employing the relationship among the roots of the cubic equation, specifically $\displaystyle\lambda_1\lambda_2\lambda_3=-\frac{16\mu^3R^2-432R^4}{4}$, and recognizing that $\mu\leq0$ leads to $\lambda_1>0$ (in accordance with the previously established conditions),
   we can conclude that for
   $4\lambda^3+\mu^2\lambda^2+72\mu R^2\lambda+16\mu^3R^2-432R^4>0$ to hold, it is necessary that
   \begin{align*}
     \lambda>\frac{-\mu^2-(Y_1^{1/3}+Y_2^{1/3})}{12}.
   \end{align*}

   {\bf Case 3.} If $\Delta_{\Delta}=0$, that is $\mu^3=-108R^2$, the equation  $4\lambda^3+\mu^2\lambda^2+72\mu R^2\lambda+16\mu^3R^2-432R^4=0$ possesses three distinct real solutions, denoted as  $\lambda_1, \lambda_2, \lambda_3$, and
   \begin{align*}
   \lambda_1=\frac{15}{2^{2/3}}R^{\frac{4}{3}},\qquad \lambda_2=\lambda_3=-3(2R)^{\frac{4}{3}}.
   \end{align*}
   In this scenario, for $4\lambda^3+\mu^2\lambda^2+72\mu R^2\lambda+16\mu^3R^2-432R^4>0$ to hold true, it is necessary that
   \begin{align*}
     \lambda>\frac{15}{2^{2/3}}R^{\frac{4}{3}}.
   \end{align*}
   Notably, this result is consistent with and encompasses the conclusions drawn from the previous cases (Case 1 and Case 2) when the condition  $\mu^3=-108R^2$ is imposed.
  \end{proof}

  \begin{figure}[htpb]
	\centering
	\setlength{\abovecaptionskip}{0.cm}
    \subfigure[]{\includegraphics[width=0.4\textwidth]{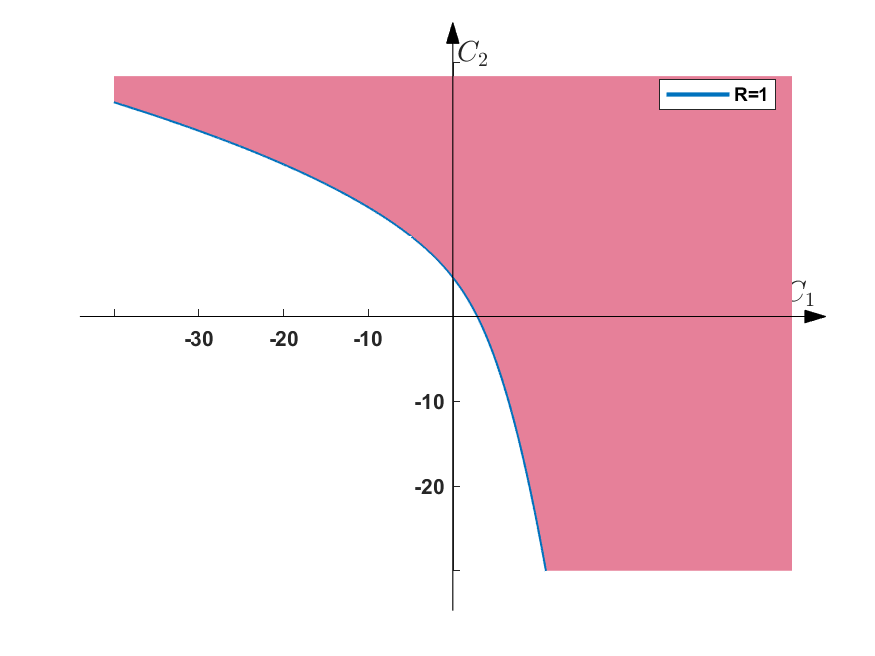}}
	\subfigure[]{\includegraphics[width=0.4\textwidth]{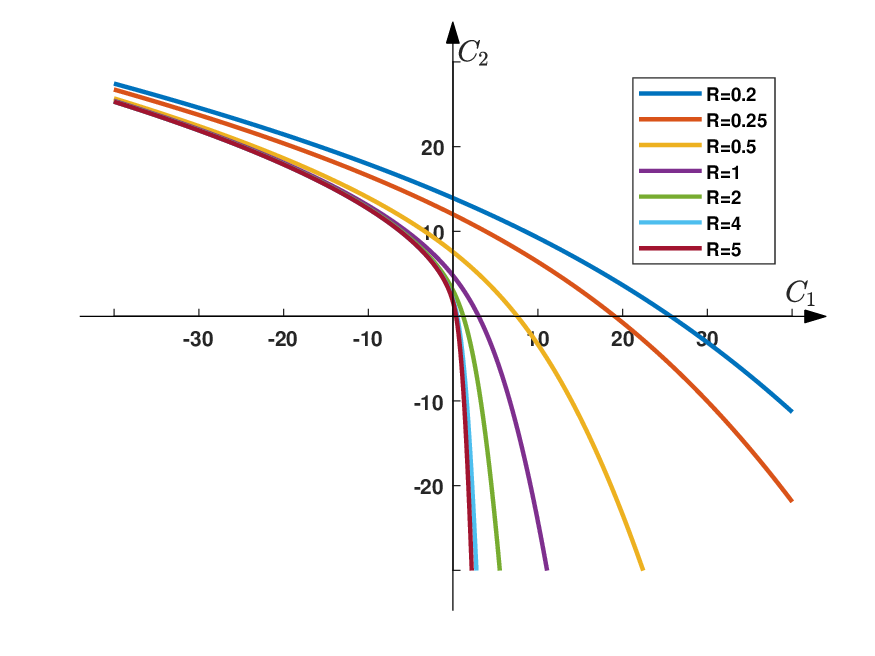}}
	\caption{The ranges for $C_1$ and $C_2$ determined by Proposition \ref{prop-bd}. }
	\label{Fig:dir}
\end{figure}
  \begin{rem}
     In Figure \ref{Fig:dir}, we illustrate the ranges  for $C_1$ and $C_2$ as specified by Proposition \ref{prop-bd}, taking into account various values of $R$. These curves depicted in the figure intersect the $C_1$ axis  at $C_1=3R^{-4/3}$ and touch the $C_2$ axis at $C_2=3(2/R)^{2/3}$. This visualization provides a clear understanding of how the ranges for $C_1$ and $C_2$ vary with different values of $R$, as prescribed by the proposition.
  \end{rem}
  \begin{thm} Suppose that $\frac{B^3}{R^2}-C_1B^2-C_2B+4=0$ admits two distinct positive roots, denoted as $A_1$ and $A_2$,  with $A_1>A_2$. Then,
  the progression angle $\Lambda^{\Theta}$ within one period of the curvature can be expressed as follows
  \begin{align*}
      2\sqrt{\frac{4C_2}{R^2}+C_1^2}\int^{A_1}_{A_2}\frac{B^{-\frac{1}{2}}}{(C_2-B^2/R^2+C_1B)\sqrt{C_2-B^2/R^2-4B^{-1}+C_1B}}dB.
  \end{align*}
  \end{thm}
  \begin{proof}
  It is straightforward to observe that within one period, the arc length $T$ is given by
  \begin{align*}
    T=2\int^{A_1}_{A_2}\frac{1}{\sqrt{C_2-B^2/R^2-4B^{-1}+C_1B}}dB
  \end{align*}
  after utilizing \eqref{BC-eq2}. Furthermore, by combining equations \eqref{BC-eq1} and \eqref{killing1}, we discover
  \begin{align*}
    \mathbf{J}=-2B^{-\frac{1}{2}}\mathbf{T}+B_s\boldsymbol{\epsilon}
  \end{align*}
  constitutes a Killing vector field along the curve $\mathbf{x}$.
  Now we can employ the spherical coordinates $\mathbf{x}(\theta,\psi)=R(\cos\psi,\cos\theta\sin\psi,\sin\theta\sin\psi)$,
  so that its equator gives the only integral geodesic of $\mathbf{J}:\mathbf{x}_{\theta}=b\mathbf{J}$, which generates
  \begin{align*}
    \mathbf{x}_{\theta}^2=R^2\sin^2\psi=b^2(4B^{-1}+B_s^2).
  \end{align*}
  From the definition of the tangent vector, we have
  \begin{align*}
      \mathbf{T}=\mathbf{x}_{\theta}\frac{d\theta}{ds}+\mathbf{x}_{\psi}\frac{d\psi}{ds}.
  \end{align*}
  Taking the inner product gives
  \begin{align*}
    \langle\mathbf{T}, \mathbf{x}_{\theta}\rangle = \langle\mathbf{x}_{\theta}, \mathbf{x}_{\theta}\rangle\frac{d\theta}{ds}=R^2\sin^2\psi\frac{d\theta}{ds}.
  \end{align*}
  Substituting the expression for $\mathbf{J}$ from earlier, and using the relationship between $\mathbf{T}$ and $\mathbf{x}_{\theta}$, we derive
  \begin{align*}
    \frac{d\theta}{ds}=\frac{\langle\mathbf{T},\mathbf{x}_{\theta}\rangle}{R^2\sin^2\psi}=\frac{-2B^{-\frac{1}{2}}}{b(4B^{-1}+B^2_s)}
    =\frac{-2B^{-\frac{1}{2}}}{b(C_2-B^2/R^2+C_1B)}
  \end{align*}
  and
  \begin{align}\label{or-A}
  \Lambda^{\Theta}=\int^T_0\frac{2\sqrt{B^{-1}}}{b(C_2-B^2/R^2+C_1B)}ds.
  \end{align}
  Let $\bar{\mathbf{x}}$ denote the  integral curve of vector field $\displaystyle \bar{\mathbf{T}}=\frac{\mathbf{J}}{\left|\mathbf{J}\right|}$. Specifically,
  \begin{align*}
     \bar{\mathbf{T}}=\frac{1}{\sqrt{C_2-B^2/R^2+C_1B}}(-2B^{-\frac{1}{2}}\mathbf{T}+B_s\boldsymbol{\epsilon}).
  \end{align*}
  Therefore at the point where $B_s=0$, we need to find the derivative
  \begin{align*}
    \bar{\mathbf{T}}_{\bar{s}}\;=&\frac{ds}{d\bar{s}}\Bigg(\frac{-2B^{-\frac{1}{2}}}{\sqrt{C_2-B^2/R^2+C_1B}}(B^{-\frac{3}{2}}\boldsymbol{\epsilon}-\frac{1}{R^2}\mathbf{x})\\
    &\qquad\qquad\qquad\qquad +\frac{B_{ss}}{\sqrt{C_2-B^2/R^2+C_1B}}\boldsymbol{\epsilon}\Bigg).
  \end{align*}
  At that given point, we have $\mathbf{x}=\bar{\mathbf{x}}$, which implies
  \begin{align*}
     \frac{ds}{d\bar{s}}=\frac{\sqrt{C_2-B^2/R^2+C_1B}}{-2B^{-\frac{1}{2}}}.
  \end{align*}
  Accordingly,
  \begin{align*}
    |\bar{\mathbf{T}}_{\bar{s}}|=\sqrt{\frac{(B_{ss}-2B^{-2})^2}{4B^{-1}}+\frac{1}{R^2}},
  \end{align*}
  which is the curvature of $\bar{\mathbf{x}}$ at the point. Since $\bar{\mathbf{x}}$ represents a planar circle,
its curvature is equivalently expressed as the reciprocal of the radius.
  Thus
  \begin{align*}
    \frac{1}{|\mathbf{x}_{\theta}|}=\sqrt{\frac{(B_{ss}-2B^{-2})^2}{4B^{-1}}+\frac{1}{R^2}}
  \end{align*}
  which generates
  \begin{align*}
    \frac{1}{b}=\sqrt{\frac{C_2}{R^2}+\frac{C_1^2}{4}}.
  \end{align*}
  Incorporating this into \eqref{or-A}  results in the anticipated  outcome.
  \end{proof}
  In light of the identity
    \begin{align*}
         \frac{1}{C_2-B^2/R^2+C_1B}\;=&\;\frac{1}{2}\left(C_2+\left(\frac{RC_1}{2}\right)^2\right)^{-1/2}\times\\
          &\quad\left(\left(\sqrt{C_2+\left(\frac{RC_1}{2}\right)^2}-\left(\frac{B}{R}-\frac{RC_1}{2}\right)\right)^{-1}\right.\\
          &\qquad\; +\left.\left(\sqrt{C_2+\left(\frac{RC_1}{2}\right)^2}+\left(\frac{B}{R}-\frac{RC_1}{2}\right)\right)^{-1}\right),
    \end{align*}
    we observe that
    \begin{eqnarray}\label{ang}
    \begin{aligned}
       \Lambda^{\Theta}=&\;\frac{2}{\sqrt{R}}\int^{\tilde{A}_1}_{\tilde{A}_2}\frac{d\tilde{B}}{(r-\tilde{B})\sqrt{(\tilde{A}_1-\tilde{B})(\tilde{B}-\tilde{A}_2)(\tilde{B}-\tilde{A}_3)}}\\
       &\qquad\quad  +\frac{2}{\sqrt{R}}\int^{\tilde{A}_1}_{\tilde{A}_2}\frac{d\tilde{B}}{(r+\tilde{B})\sqrt{(\tilde{A}_1-\tilde{B})(\tilde{B}-\tilde{A}_2)(\tilde{B}-\tilde{A}_3)}},
    \end{aligned}
    \end{eqnarray}
  where $\displaystyle r=\sqrt{C_2+\left(\frac{RC_1}{2}\right)^2}$, $\displaystyle \tilde{B}=\frac{B}{R}-\frac{RC_1}{2}$, and $\displaystyle \tilde{A}_i=\frac{A_i}{R}-\frac{RC_1}{2}$ $(i=1, 2, 3)$.
  Additionally, the expressions of $A_1, A_2$ and $A_3$ can be derived as follows
  \begin{align*}
     \theta&=\;\arccos\frac{108-9C_1C_2R^2-2C_1^3R^4}{2R(C_1^2R^2+3C_2)^{3/2}},\\
     A_1&=\;\frac{C_1R^2}{3}+\frac{2R}{3}\sqrt{C_1^2R^2+3C_2}\cos\frac{\theta-\pi}{3},\\
     A_2&=\;\frac{C_1R^2}{3}+\frac{2R}{3}\sqrt{C_1^2R^2+3C_2}\cos\frac{\theta+\pi}{3},\\
     A_3&=\;\frac{C_1R^2}{3}-\frac{2R}{3}\sqrt{C_1^2R^2+3C_2}\cos\frac{\theta}{3}.
  \end{align*}
  Let $D$ denote the lower bounds of $\lambda$ as specified in Proposition \ref{prop-bd}. Notably, $\theta$ exhibits a monotonic increase with respect to $\lambda$, and its limiting behavior is characterized by
  \begin{equation}\label{theta-lim}
  \lim_{\lambda\to D}\theta=0,\qquad \lim_{\lambda\to+\infty}\theta=\frac{\pi}{2}.
  \end{equation}
  For the sake of brevity and clarity, we introduce the notations
  \begin{align*}
     d&=\;\sqrt{C_1^2R^2+3C_2},\\
     a&=\;\frac{2}{3}d\cos\frac{\theta-\pi}{3}-\frac{RC_1}{6},\\
     b&=\;\frac{2}{3}d\cos\frac{\theta+\pi}{3}-\frac{RC_1}{6},\\
     c&=\;-\frac{2}{3}d\cos\frac{\theta}{3}-\frac{RC_1}{6},
  \end{align*}
  which exhibit the subsequent characteristics.
  \begin{lem}\label{lem-inq}  $a,r,d$ is monotonically increasing and $b,c$ is monotonically decreasing with respect to $C_2$. moreover, $a>b>c$ and $r>b$.
  \end{lem}
  \begin{lem}\label{lem-lim}
   The limits stated below are valid.
  \begin{align*}
     \lim_{\lambda\to+\infty}\frac{a}{r}=1,\quad \lim_{\lambda\to+\infty}\frac{c}{a}=-1, \quad \lim_{\lambda\to+\infty}\frac{b}{a}=0, \quad \lim_{\lambda\to+\infty}a^2(r-a)R=2.
  \end{align*}
  \end{lem}
Upon performing a series expansion centered at $x=0$ with the parameter $x=1/(d+1)$, we establish
\begin{lem}
 $\displaystyle \frac{(a-b)(c-r)}{(a-c)(b-r)}$ is monotonically increasing with respect to $C_2$, and furthermore, $\displaystyle \frac{(a-b)(c-r)}{(a-c)(b-r)}<1$.
\end{lem}
  Based on \cite{gr}, we obtain an elegant expression for the given integral
   \begin{align*}
       \int^{a}_{b}\frac{d\tilde{B}}{(r+\tilde{B})\sqrt{(a-\tilde{B})(\tilde{B}-b)(\tilde{B}-c)}}=\frac{2}{(a+r)\sqrt{a-c}}\Pi\left(\frac{a-b}{a+r},\frac{\pi}{2},\frac{a-b}{a-c}\right),
    \end{align*}
    where $\Pi$ denotes the elliptic integral of the third kind defined as follows
    \begin{align*}
     \Pi(n,\varphi,m)&=\int^{\varphi}_0\frac{d\theta}{(1-n\sin^2\theta)\sqrt{1-m\sin^2\theta}}\\
     &=\int_0^{\sin\varphi}\frac{dt}{(1-nt^2)\sqrt{(1-mt^2)(1-t^2)}}.
    \end{align*}
    Utilizing  Lemma \ref{lem-inq}, we confirm that $r>b$ holds for all $R, C_1$ and $C_2$. Further, drawing from  \cite{gr}, we express the given integral as
    \begin{align*}
      \int^{a}_{b}\frac{d\tilde{B}}{(r-\tilde{B})\sqrt{(a-\tilde{B})(\tilde{B}-b)(\tilde{B}-c)}}\;=&
      \frac{-2(c-b)}{(c-r)(b-r)\sqrt{a-c}}\times\\
      &\;\Pi\left(\frac{(a-b)(c-r)}{(a-c)(b-r)},\frac{\pi}{2},\frac{a-b}{a-c}\right)\\
                   &+\frac{-2}{(c-r)\sqrt{a-c}}\Pi\left(0,\frac{\pi}{2},\frac{a-b}{a-c}\right).
    \end{align*}
 Subsequently, we derive the refined form of $\Lambda^{\Theta}$,
 \begin{equation}
 \begin{aligned}
   \Lambda^{\Theta}\;=\;&\frac{4}{(a+r)\sqrt{R(a-c)}}\Pi\left(\frac{a-b}{a+r},\frac{\pi}{2},\frac{a-b}{a-c}\right)\\
                    &+\frac{4(b-c)}{(r-c)(r-b)\sqrt{R(a-c)}}\Pi\left(\frac{(a-b)(c-r)}{(a-c)(b-r)},\frac{\pi}{2},\frac{a-b}{a-c}\right)\\
                   &+\frac{4}{(r-c)\sqrt{R(a-c)}}\Pi\left(0,\frac{\pi}{2},\frac{a-b}{a-c}\right). \label{anP}
 \end{aligned}
 \end{equation}
 \begin{prop}\label{prop-ub}
 As $\lambda\to D$, the limit $\Lambda^{\Theta}$ of is given by
    \begin{align*}
      \lim_{\lambda\to D}\Lambda^{\Theta}\;=&\;\frac{12\pi}{\sqrt{Rd_0}}\Bigg(\frac{1}{2d_0+6r_0-RC_1}+\frac{1}{6r_0+4d_0+RC_1}\\
      &\qquad\qquad\qquad +\frac{6d_0}{(6r_0+4d_0+RC_1)(6r_0-2d_0+RC_1)}\Bigg),
    \end{align*}
    where $\displaystyle r_0=\sqrt{\frac{D}{R^2}+\frac{R^2C_1^2}{4}}$ and $d_0=\sqrt{\frac{3D}{R^2}+C_1^2R^2}$, with $D$ being the lower bound of $\lambda$ in Proposition \ref{prop-bd}.
 \end{prop}
 \begin{proof}
  Observe that
   \begin{gather*}
     \lim_{\lambda\to D}r=r_0,\quad \lim_{\lambda\to D}d=d_0, \quad\lim_{\lambda\to D}c=-\frac{2d_0}{3}-\frac{RC_1}{6}, \\
     \lim_{\lambda\to D}a=\lim_{\lambda\to D}b=\frac{d_0}{3}-\frac{RC_1}{6}.
   \end{gather*}
  Substituting these limits into the expression \eqref{anP} and simplifying the resulting terms, we directly obtain the stated result.
 \end{proof}
 \begin{prop}
 The following limit holds true:
     \begin{align*}
      \lim_{\lambda\to +\infty}\Lambda^{\Theta}\;=\pi.
    \end{align*}
 \end{prop}
 \begin{proof}
 According to Lemma \ref{lem-lim}, it is immediately evident that
 \begin{align*}
   \lim_{\lambda\to+\infty}\frac{4}{(a+r)\sqrt{R(a-c)}}\Pi\left(\frac{a-b}{a+r},\frac{\pi}{2},\frac{a-b}{a-c}\right)=0
 \end{align*}
 and similarly,
 \begin{align*}
   \lim_{\lambda\to+\infty}\frac{4}{(r-c)\sqrt{R(a-c)}}\Pi\left(0,\frac{\pi}{2},\frac{a-b}{a-c}\right)=0.
 \end{align*}
 Now, let us focus our attention to the remaining component within \eqref{anP}. Supposing that
 \begin{align*}
   \alpha^2=\frac{(a-b)(c-r)}{(a-c)(b-r)},\quad k^2=\frac{a-b}{a-c},
 \end{align*}
we may readily observe that
 \begin{align*}
   \lim_{\lambda\to+\infty}\alpha^2=1,\qquad \lim_{\lambda\to+\infty}k^2=\frac{1}{2}.
 \end{align*}
 According to \cite{bf}, we obtain the expansion
 \begin{align*}
   \int^{\frac{\pi}{2}}_0\frac{d\theta}{(1-\alpha^2\sin^2\theta)\sqrt{1-k^2\sin^2\theta}}\;=&\;\int^{\frac{\pi}{2}}_0\frac{d\theta}{\sqrt{1-k^2\sin^2\theta}}\\
   &\quad\qquad+\frac{1}{1-k^2}\int^{\frac{\pi}{2}}_0\sqrt{1-k^2\sin^2\theta}d\theta\\
   &\quad\qquad+\frac{(2-k^2(1+\alpha^2))\pi}{4(1-k^2)\sqrt{1-k^2}\sqrt{1-\alpha^2}}\\
   &\quad\qquad +o(1-\alpha^2),
 \end{align*}
 where $\alpha^2\neq1$ and $o(1-\alpha^2)$ signifies that the remaining terms are bounded by $(1-\alpha^2)C$ for some constant $C$.
 One can verify that
 \begin{align*}
 \lim_{k^2\to\frac{1}{2}}\int^{\frac{\pi}{2}}_0\frac{d\theta}{\sqrt{1-k^2\sin^2\theta}}=0, \quad
 \lim_{k^2\to\frac{1}{2}}\int^{\frac{\pi}{2}}_0\sqrt{1-k^2\sin^2\theta}d\theta=0.
 \end{align*}
  Direct computation shows
  \begin{align*}
     \lim_{\lambda\to+\infty}\frac{4(b-c)}{(r-c)(r-b)\sqrt{R(a-c)}\sqrt{1-\alpha^2}}=\sqrt{2},
  \end{align*}
  and
  \begin{align*}
   \lim_{\lambda\to+\infty}\frac{(2-k^2(1+\alpha^2))}{4(1-k^2)\sqrt{1-k^2}}=\frac{\sqrt{2}}{2}.
  \end{align*}
  Given these limits, it is direct to deduce the desired result.
 \end{proof}

 \begin{figure}[htpb]
\centering
\setlength{\abovecaptionskip}{0.cm}
\subfigure[]{\includegraphics[width=0.3\textwidth]{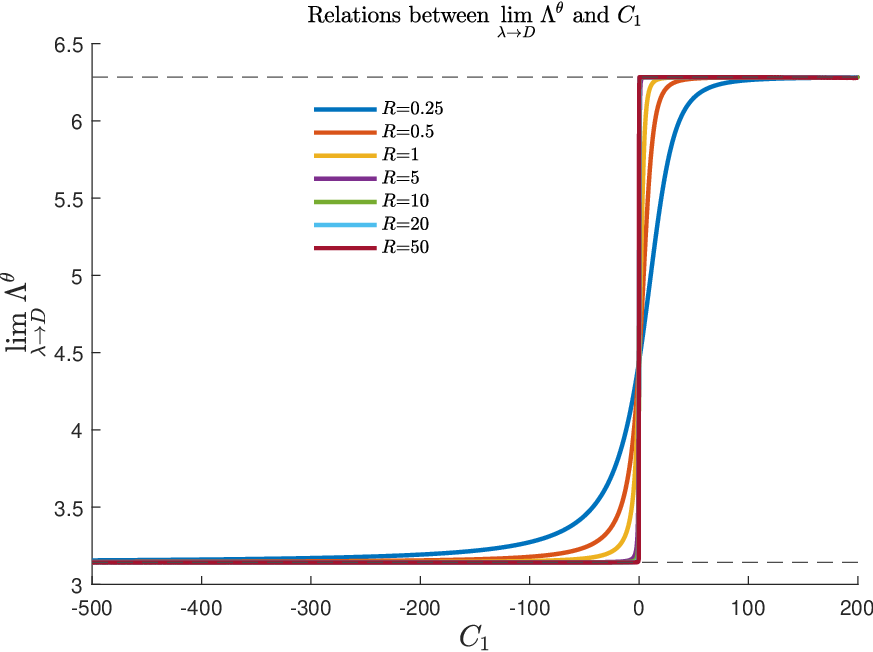}}\qquad
\subfigure[]{\includegraphics[width=0.315\textwidth]{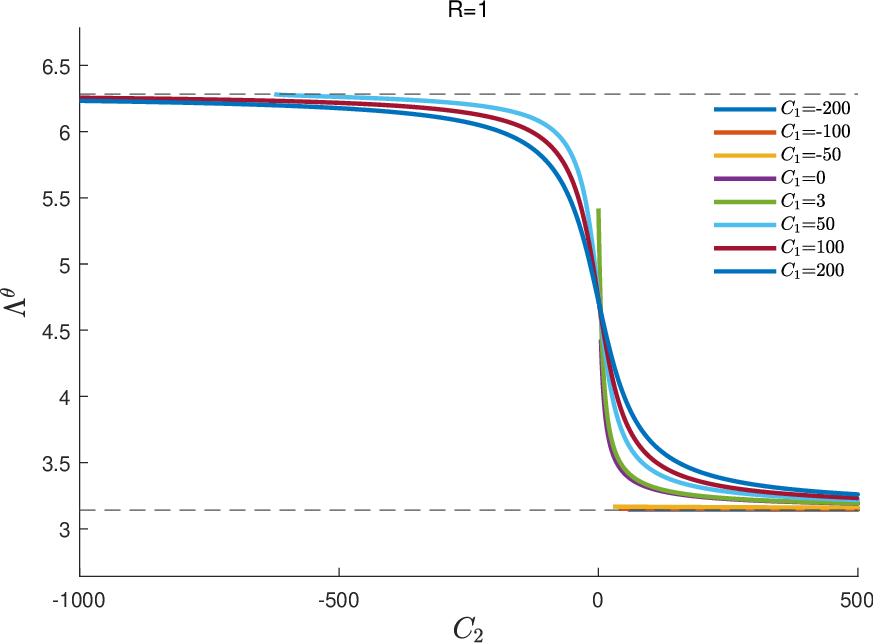}}
\caption{The incremental or decremental processes  of $\Lambda^{\Theta}$.}
\label{Fig:trend}
\end{figure}
Fig. \ref{Fig:trend} presents a comprehensive visualization showcasing the varying trends of the progression angle $\displaystyle \Lambda^{\Theta}$, contingent upon diverse variables and inputs.
Specifically, Fig. \ref{Fig:trend} (a), by modulating the parameter $R$, it becomes evident that $\displaystyle \lim_{\lambda\to D}\Lambda^{\Theta}$ exhibits a monotonic increase with respect to $C_1$.
Additionally, two notable limits can be discerned: $\displaystyle \lim_{C_1\to-\infty}\lim_{\lambda\to D}\Lambda^{\Theta}=\pi$ and  $\displaystyle \lim_{C_1\to+\infty}\lim_{\lambda\to D}\Lambda^{\Theta}=2\pi$.

Utilizing \eqref{BC-eq2} for analysis, upon investigation of the consequences of resizing the sphere, we ascertain that the scale of the spherical radius $R$ does not influence the fundamental dynamics of the progression angle $\Lambda^{\Theta}$.
As depicted in  Fig. \ref{Fig:trend} (b),
 when $C_1$ is held at different constants, the progression angle $\Lambda^{\Theta}$ consistently demonstrates a monotonic decrease in response to variations in $C_2$.
To substantiate this observation, in Appendix B, we provide a rigorous proof of the monotonicity of the progression angle $\Lambda^{\Theta}$ specifically for the case where $C_1=0$ and $R=1$,
 serving as an illustrative example.
 Analogous approaches can be adopted to delve into the intricate behaviors of $\Lambda^{\Theta}$ under diverse parameter settings and conditions,
 offering a comprehensive understanding of its dynamic properties across various scenarios.
 Indeed, we have firmly  established the following result.
 \begin{prop}
   \begin{itemize}
      \item[(1)]  For $C_1=0$, we have  $\displaystyle \lim_{\lambda\to D}\Lambda^{\Theta}=\sqrt{2}\pi.$
      \item[(2)]  For all $R$, $\displaystyle\lim_{\lambda\to D}\Lambda^{\Theta}$ is monotonic increasing with respect to $C_1$, and
      \begin{align*}
        \lim_{C_1\to-\infty}\lim_{\lambda\to D}\Lambda^{\Theta}=\pi,\quad \lim_{C_1\to+\infty}\lim_{\lambda\to D}\Lambda^{\Theta}=2\pi.
      \end{align*}
      \item[(3)] Holding $C_1$ and $R$ constant, $\Lambda^{\Theta}$ is monotonic decreasing with respect to $C_2$ and
      $\displaystyle
        \pi<\Lambda^{\Theta}<\lim_{\lambda\to D}\Lambda^{\Theta}.
      $
       \item[(4)] For all $R, C_1, C_2$, the range of $\Lambda^{\Theta}$ is constrained to $\pi<\Lambda^{\Theta}<2\pi$.
    \end{itemize}
 \end{prop}
\begin{thm}\label{thm-ca}
Let $\mathbf{x}$  be a smoothly embedded closed curve residing on a sphere, whose curvatures satisfy the condition stated in \eqref{BC-eq1}. Then $\mathbf{x}$  is inherently a planar curve. In other words, the smoothly closed and embedded
generalized equi-centro-affine extremal curves residing on a sphere are precisely planar circles.
\end{thm}
Regarding the classification of closed generalized equi-centro-affine extremal curves on a sphere, we have
\begin{thm}
Let $\mathbf{x}$ be a closed generalized equi-centro-affine extremal curve on sphere. Then we have the following possibilities for $\mathbf{x}$.
 \begin{itemize}
   \item[(1)] $\mathbf{x}$ is a planar circle on sphere;
   \item[(2)]$\mathbf{x}=\mathbf{x}_{p,q}$, a curve whose curvature is a noncanstant perodic solution of \eqref{BC-eq1} depending upon $(p,q)\in \mathbf{N}\times \mathbf{N}$. The pair $(p,q)$ is by no means arbitrary and must be such that $p/q$ in defined to the open interval $\displaystyle\left(1/2,\; \lim_{\lambda\to D}\Lambda^{\Theta}/2\pi\right)$, where $\displaystyle\lim_{\lambda\to D}\Lambda^{\Theta}$ is obtained in Proposition \ref{prop-ub}.
 \end{itemize}
\end{thm}

 \begin{figure}[htpb]
\centering
\setlength{\abovecaptionskip}{0.cm}
\subfigure[$p=2,\;q=3,\;C_1=0$]{\includegraphics[width=0.3\textwidth]{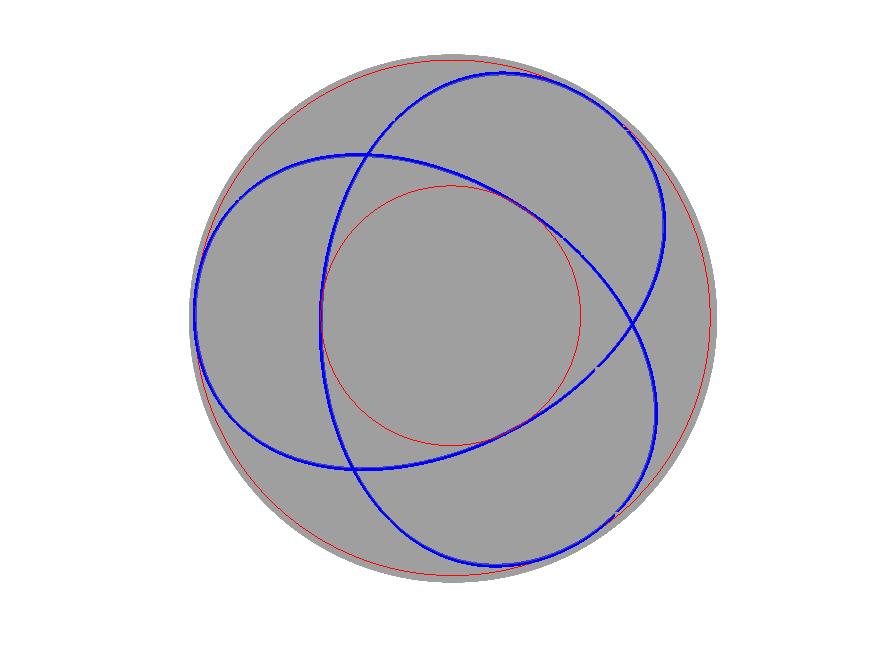}}
\subfigure[$p=3,\;q=5,\;C_1=0$]{\includegraphics[width=0.3\textwidth]{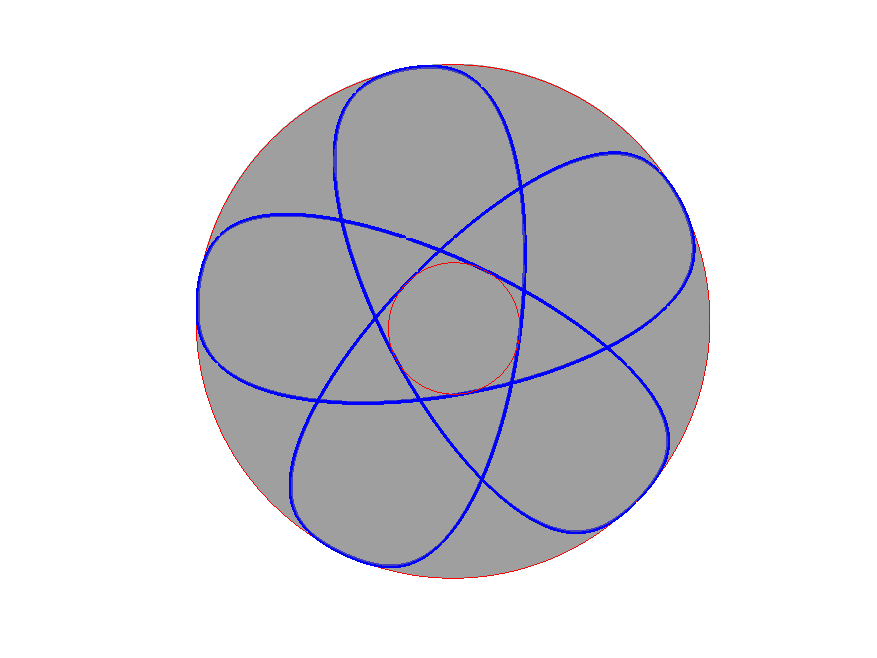}}
\subfigure[$p=4,\;q=7,\;C_1=0$]{\includegraphics[width=0.3\textwidth]{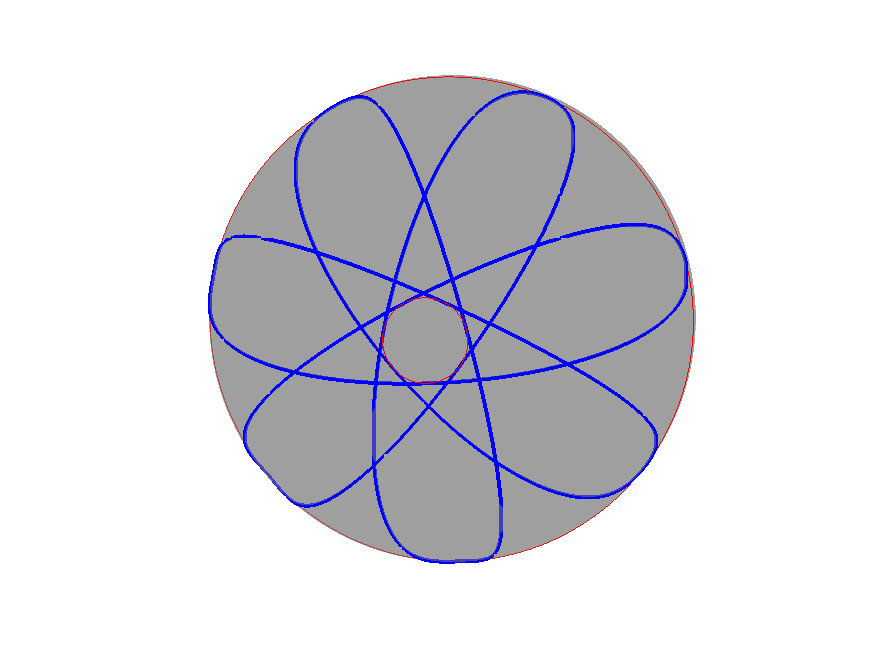}}
\subfigure[$p=3,\;q=4,\;C_1=20$]{\includegraphics[width=0.3\textwidth]{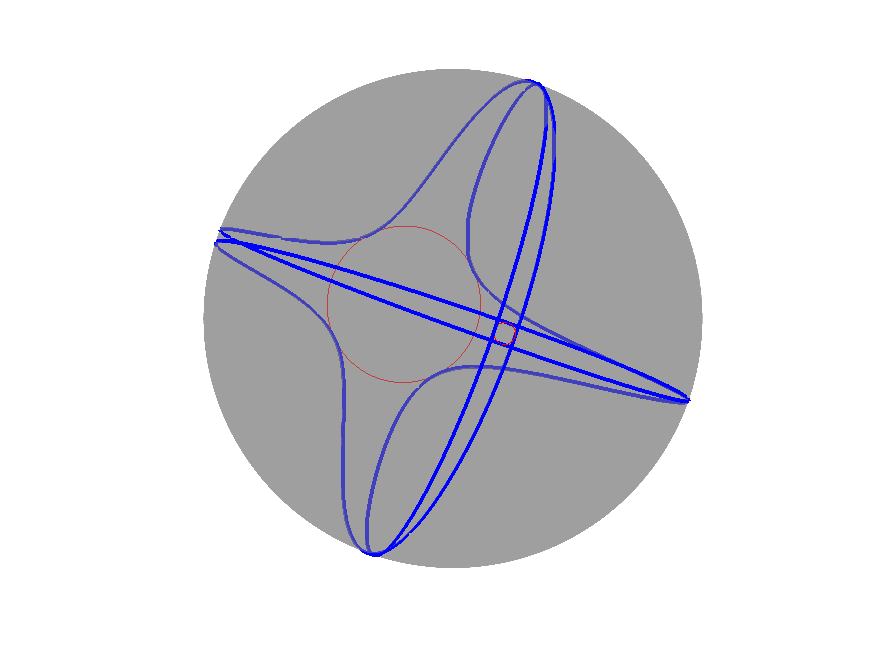}}
\subfigure[$p=4,\;q=5,\;C_1=30$]{\includegraphics[width=0.33\textwidth]{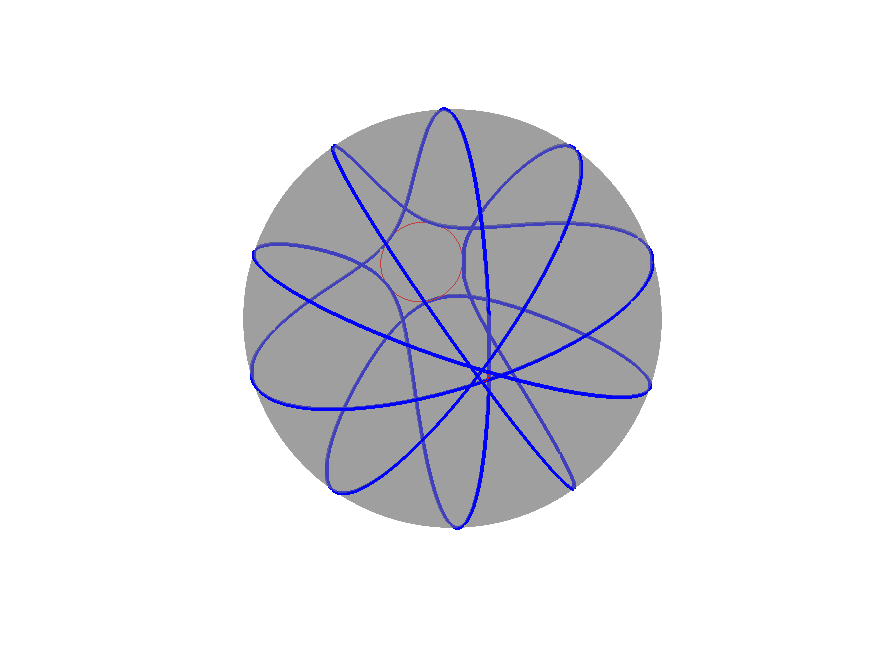}}
\subfigure[$p=5,\;q=6,\;C_1=30$]{\includegraphics[width=0.32\textwidth]{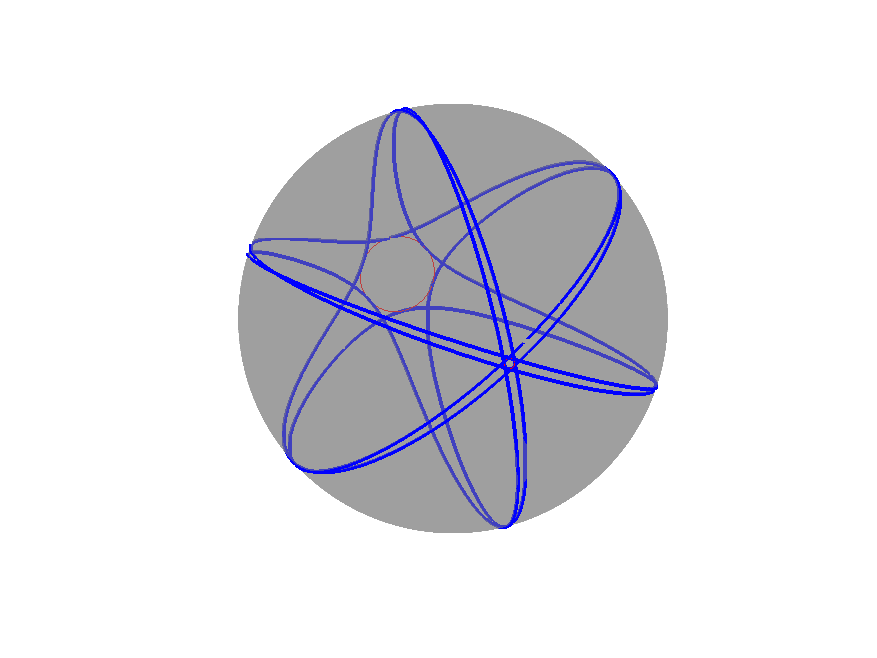}}
\caption{The closed extremal curves on unit sphere.}
\label{Fig:closed}
\end{figure}
Taking $R=1$ and $C_1=0$  as an illustrative example, we demonstrate the classification of closed generalized equi-centro-affine extremal curves on a sphere, which restricts the ratio $p/q$ to lie within the open interval $(1/2, \sqrt{2}/2)$.
To extend this range and obtain closed curves with $p/q$ values in the interval $(\sqrt{2}/2,\;1)$, we can adjust the parameter $C_1$ to specific values such as $C_1=20$ or $C_1=30$.
As depicted in Fig. \ref{Fig:closed}, several closed curves $\mathbf{x}_{p,q}$ are presented, all of which can be considered as generalized equi-centro-affine extremal curves on the sphere under these parameter settings.

\appendix
{\centering\section*{Appendix A: some calculations for Theorem \ref{thm-second}}}
\setcounter{equation}{0}
\setcounter{subsection}{0}
\setcounter{thm}{0}
\renewcommand{\theequation}{A.\arabic{equation}}
\renewcommand{\thesubsection}{A.\arabic{subsection}}
\renewcommand{\thethm}{A.\thesection\arabic{thm}}
In this part, we provide a more detailed exposition of the computation process leading to the derivation of the second variational formula presented in Theorem \ref{thm-second}.
By \eqref{srt}, we have
\begin{align*}
 \frac{\partial}{\partial t}&\left(S_n^{-\frac{n+1}{n+2}}(\frac{1}{\varrho^2}S_{n-1}-(n+1)S_1S_n)\right)\\
 =&U_{,ij}\Bigg(-\frac{n+1}{n+2}S_n^{-\frac{2n+3}{n+2}}T^{ij}_{n-1}(\frac{1}{\varrho^2}S_{n-1}-(n+1)S_1S_n)\\
 &\qquad+S_n^{-\frac{n+1}{n+2}}\left(\frac{1}{\varrho^2}T^{ij}_{n-2}-(n+1)(T^{ij}_0S_n+S_1T^{ij}_{n-1})\right)\Bigg)\\
 &\qquad\quad+U\Bigg(-\frac{n+1}{n+2}S_n^{-\frac{2n+3}{n+2}}\left(\frac{1}{\varrho^4}S_{n-1}^2-(n+1)S_1^2S_{n}^2-\frac{n}{\varrho^2}S_1S_{n-1}S_n\right)\\
       &\quad\qquad\quad+S_n^{-\frac{n+1}{n+2}}\bigg(\frac{2}{\varrho^4}S_{n-2}-\frac{n}{\varrho^2}S_n+2(n+1)S_2S_n\\
       &  \quad\qquad\qquad           -\frac{n(n+1)}{\varrho^2}S_n-2(n+1)S_1^2S_n-\frac{n}{\varrho^2}S_1S_{n-1}\bigg)\Bigg)\\
 &\qquad\qquad\qquad+W^j\Bigg(-\frac{n+1}{n+2}S_n^{-\frac{2n+3}{n+2}}S_{n,j}\left(\frac{1}{\varrho^2}S_{n-1}-(n+1)S_1S_n\right)\\
 &\qquad\qquad\qquad\qquad+S^{-\frac{n+1}{n+2}}\left(\frac{1}{\varrho^2}S_{n-1,j}-(n+1)(S_{1,j}S_n+S_1S_{n,j})\right)\Bigg).
\end{align*}
On the other hand,
\begin{align*}
\frac{\partial}{\partial t}\left(\left(S_n^{-\frac{n+1}{n+2}}\right)_{,ij}T^{ij}_{n-1}\right)=\frac{\partial}{\partial t}\left(\left(S_n^{-\frac{n+1}{n+2}}\right)_{,ij}\right)T^{ij}_{n-1}+\left(S_n^{-\frac{n+1}{n+2}}\right)_{,ij}\frac{\partial}{\partial t}T^{ij}_{n-1},
\end{align*}
and
\begin{align*}
  \frac{\partial}{\partial t}\left(\left(S_n^{-\frac{n+1}{n+2}}\right)_{,ij}\right)&=\bar{\nabla}_{\frac{\partial}{\partial t}}\bar{\nabla}_{e_j}e_i\left(S_n^{-\frac{n+1}{n+2}}\right)-\frac{\partial}{\partial t}\left(\Gamma^k_{ij}e_k\left(S_n^{-\frac{n+1}{n+2}}\right)\right)\\
  &=\left(\frac{\partial}{\partial t}\left(S_n^{-\frac{n+1}{n+2}}\right)\right)_{,ij}-\frac{\partial \Gamma^k_{ij}}{\partial t}e_k\left(S_n^{-\frac{n+1}{n+2}}\right)\\
  &\qquad\qquad+\bar{R}\left(e_j,\frac{\partial \mathbf{x}}{\partial t}\right)e_i\left(S_n^{-\frac{n+1}{n+2}}\right).
\end{align*}
According to \eqref{srt}, we may find
\begin{align*}
 \left(\frac{\partial}{\partial t}S_n^{-\frac{n+1}{n+2}}\right)_{,ij}T^{ij}_{n-1}=&A_4U_{,mkij}+A_3U_{,mki}+A_2U_{,mk}+A_1U_{,m}\\
 &\quad+A_0U+B_2W^m_{,ij}+B_1W^m_{,i}+B_0W^m,
\end{align*}
where
\begin{align*}
  A_4&=-\frac{n+1}{n+2}S_n^{-\frac{2n+3}{n+2}}T^{mk}_{n-1}T^{ij}_{n-1},\\
  A_3&=-\frac{2(n+1)}{n+2}\left(S_n^{-\frac{2n+3}{n+1}}T^{mk}_{n-1}\right)_{,j}T^{ij}_{n-1},\\
  A_2&=-\frac{n+1}{n+2}\left(S_n^{-\frac{2n+3}{n+2}}T^{mk}_{n-1}\right)_{,ij}T^{ij}_{n-1}-\frac{n+2}{n+2}S_n^{-\frac{2n+3}{n+2}}(S_1S_n+\frac{1}{\varrho^2}S_{n-1})T^{mk}_{n-1},\\
  A_1&=-\frac{2(n+1)}{n+2}\left(S_n^{-\frac{2n+3}{n+2}}(S_1S_n+\frac{1}{\varrho^2}S_{n-1})\right)_{k}T^{mk}_{n-1},\\
  A_0&=-\frac{n+1}{n+2}\left(S_n^{-\frac{2n+3}{n+2}}(S_1S_n+\frac{1}{\varrho^2}S_{n-1})\right)_{,ij}T^{ij}_{n-1},\\
  B_2&=\left(S_n^{-\frac{n+1}{n+2}}\right)_mT^{ij}_{n-1},\\
  B_1&=2\left(S_n^{-\frac{n+1}{n+2}}\right)_{,mj}T^{ij}_{n-1},\\
  B_0&=\left(S_n^{-\frac{n+1}{n+2}}\right)_{,mij}T^{ij}_{n-1}.
\end{align*}
\eqref{Gt}, \eqref{Tt} and above computations yield Theorem  \ref{thm-second}.

\appendix
{\centering\section*{Appendix B: monotonicity of the progression angle}}

\setcounter{equation}{0}
\setcounter{subsection}{0}
\setcounter{thm}{0}
\renewcommand{\theequation}{B.\arabic{equation}}
\renewcommand{\thesubsection}{B.\arabic{subsection}}
\renewcommand{\thethm}{B.\thesection\arabic{thm}}

According to  \eqref{BC-eq2},  upon examining the effects of scaling the sphere, we ascertain that the magnitude of the spherical radius $R$ has no impact on the overarching dynamics of the progression angle $\Lambda^{\Theta}$.
Consequently, a simplifying assumption of $R=1$ can be safely made without affecting the analysis.
Here, we present a meticulous proof of the monotonicity of the progression angle for the specific case where $R=1$ and $C_1=0$, which serves as a prototypical example.
It is worth noting that when $C_1\neq0$, the progression angle can be expanded as a series in terms of the parameter $x=1/(d+1)$ around $x=0$, facilitating further analysis in more complex scenarios.
For $R=1, C_1=0$, we have $r>a$,
and \eqref{ang} can be rewritten as
\begin{equation}
\begin{aligned}
   \Lambda^{\Theta}\;=\;&\frac{4}{(a+r)\sqrt{(a-c)}}\Pi\left(\frac{a-b}{a+r},\frac{\pi}{2},\frac{a-b}{a-c}\right)\\
                    &-\frac{4}{(a-r)\sqrt{(a-c)}}\Pi\left(\frac{a-b}{a-r},\frac{\pi}{2},\frac{a-b}{a-c}\right).\label{eint}
 \end{aligned}
 \end{equation}
Utilizing the series expansions for the elliptic integral of the third kind, as detailed in the reference \cite{bf}, under the conditions where $\alpha<-1,\; 0\leq k<1$ or $0<\alpha<1,\; 0\leq k<\alpha$,
we can express the elliptic integral as follows:
\begin{align}\label{Ser1}
  \Pi\left(\alpha,\frac{\pi}{2},k\right)=\sum_{m=0}^{\infty}c_m k^m,
\end{align}
 where
 \begin{eqnarray*}
 \begin{aligned}
    c_0=&\;\frac{\pi}{2\sqrt{1-\alpha}}, \quad c_1=\frac{\pi}{4\alpha}\left(\frac{1}{\sqrt{1-\alpha}}-2\right),\\
    c_2=&\;\frac{3\pi}{32\alpha^2}\left(\frac{2}{\sqrt{1-\alpha}}-2-\alpha\right),\\
    c_3=&\;\frac{5\pi}{256\alpha^3}\left(-4\alpha-3\alpha^2-8+\frac{8}{\sqrt{1-\alpha}}\right),\\
    2(m+1)\alpha c_{m+1}=&\;\frac{\pi}{2(2m-1)}
    \left(
     \begin{array}{c}
       -\frac{1}{2} \\
       m \\
     \end{array}
    \right)^2+(1-2m)c_{m-1}+(2m+1+2m\alpha)c_m.
 \end{aligned}
 \end{eqnarray*}

For the first part of $\Lambda^{\Theta}$, i.e., $\displaystyle \frac{4}{(a+r)\sqrt{(a-c)}}\Pi\left(\frac{a-b}{a+r},\frac{\pi}{2},\frac{a-b}{a-c}\right)$, we derive
\begin{align*}
   \frac{4c_0}{(a+r)\sqrt{a-c}}\;&=\;\pi r^{-\frac{3}{2}}-\frac{3\pi}{2}r^{-\frac{9}{2}}+\frac{79\pi}{8}r^{-\frac{15}{2}}+O(r^{-\frac{21}{2}}),\\
   \frac{4c_1}{(a+r)\sqrt{a-c}}\left(\frac{a-b}{a-c}\right)\;&=\;\frac{(2-\sqrt{2})\pi}{4}r^{-\frac{3}{2}}-\frac{5\pi}{4}r^{-\frac{9}{2}}\\
                         &\qquad\qquad+\frac{115-36\sqrt{9}}{16}\pi r^{-\frac{15}{2}}+O(r^{-\frac{21}{2}}),\\
   \frac{4c_2}{(a+r)\sqrt{a-c}}\left(\frac{a-b}{a-c}\right)^2\;&=\;\frac{(24-15\sqrt{2})\pi}{64}r^{-\frac{3}{2}}-\frac{(42-15\sqrt{2})\pi}{32}r^{-\frac{9}{2}}\\
                         &\qquad\qquad+\frac{477-171\sqrt{2}}{64}\pi r^{-\frac{15}{2}}+O(r^{-\frac{21}{2}}),\\
   \frac{4c_3}{(a+r)\sqrt{a-c}}\left(\frac{a-b}{a-c}\right)^3\;&=\;\frac{5(64-43\sqrt{2})\pi}{1024}r^{-\frac{3}{2}}-\frac{5(72-39\sqrt{2})\pi}{256}r^{-\frac{9}{2}}\\
                         &\qquad\qquad+\frac{5(1688-791\sqrt{2})}{1024}\pi r^{-\frac{15}{2}}+O(r^{-\frac{21}{2}}).
 \end{align*}

For the second part of $\Lambda^{\Theta}$, that is, $\displaystyle \frac{-4}{(a-r)\sqrt{(a-c)}}\Pi\left(\frac{a-b}{a-r},\frac{\pi}{2},\frac{a-b}{a-c}\right),$ we see
 \begin{align*}
   \frac{-4c_0}{(a-r)\sqrt{a-c}}\;&=\;\pi+\frac{\pi}{2}r^{-3}+\frac{11\pi}{8}r^{-6}+O(r^{-9}),\\
   \frac{-4c_1}{(a-r)\sqrt{a-c}}\left(\frac{a-b}{a-c}\right)\;&=\;\frac{\sqrt{2}\pi}{4}r^{-\frac{3}{2}}-\frac{\pi}{2}r^{-3}-\frac{7\pi}{4}r^{-6}\\
                         &\qquad\qquad+\frac{9\sqrt{2}}{4}\pi r^{-\frac{15}{2}}+O(r^{-9}),\\
   \frac{-4c_2}{(a-r)\sqrt{a-c}}\left(\frac{a-b}{a-c}\right)^2\;&=\;\frac{3\sqrt{2}\pi}{64}r^{-\frac{3}{2}}-\frac{15\sqrt{2}\pi}{32}r^{-\frac{9}{2}}+\frac{3\pi}{8}r^{-6}\\
                      &\qquad\qquad-\frac{9\sqrt{2}}{64}\pi r^{-\frac{15}{2}}+O(r^{-9}),\\
   \frac{-4c_3}{(a-r)\sqrt{a-c}}\left(\frac{a-b}{a-c}\right)^3\;&=\;\frac{15\sqrt{2}\pi}{1024}r^{-\frac{3}{2}}-\frac{55\sqrt{2}\pi}{256}r^{-\frac{9}{2}}\\
                      &\qquad\qquad-\frac{955\sqrt{2}}{1024}\pi r^{-\frac{15}{2}}+O(r^{-9}).
 \end{align*}

Then
\begin{equation}
\begin{aligned}\label{app-r}
  \Lambda^{\Theta}\;=&\;\pi\left(1+\frac{280-49\sqrt{2}}{128}r^{-\frac{3}{2}}-\frac{350-35\sqrt{2}}{64} r^{-\frac{9}{2}}\right.\\
  &\qquad\qquad\left.+\frac{4193-735\sqrt{2}}{128}r^{-\frac{15}{2}}\right)+O(r^{-9}).
\end{aligned}
\end{equation}
Then by \eqref{app-r} we can easily verify that $\Lambda^{\Theta}$ is monotonically decreasing when $\displaystyle r>\left(6\sqrt{3}\right)^{1/3}$, and converges to $\pi$ as $r\to+\infty$.
\begin{figure}[htpb]
\centering
\setlength{\abovecaptionskip}{0.cm}
\setlength{\abovecaptionskip}{0.cm}
	\includegraphics[width=0.7\textwidth]{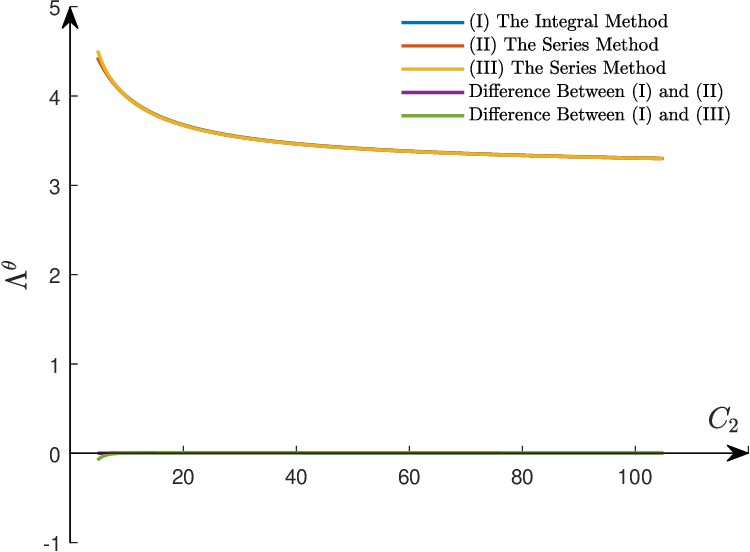}
	\caption{The comparison of different method to calculate $\Lambda^{\Theta}$.}
\label{Fig:compare}
\end{figure}

\begin{rem}
In Fig. \ref{Fig:compare},  the curves depicting $\Lambda^{\Theta}(C_2)$ are precisely rendered using three distinct methods: the numerical integral approach defined in \eqref{eint}, the summation of the initial four terms from the series expansion in \eqref{Ser1}, and the approximate series expansions \eqref{app-r} with respect to $r$. Notably, these three representations exhibit remarkable congruence, with minimal deviations among them, essentially indicating a near-perfect overlap.
 \end{rem}

\noindent {\bf Acknowledgments.}
 This work was supported by the National NSF of China Grants-12431008.

\end{document}